\documentclass[11pt]{amsart}
\usepackage{amsfonts}
\usepackage{amssymb,latexsym,amsmath,amsthm, mathrsfs}
\usepackage[all]{xy}
\usepackage{bbm}

\def\mycaption#1#2{
    \addtocounter{figure}{1}
    \parbox[b]{#1}{{\scshape Figure $\arabic{figure}$.}\ #2}
}   

\def\pd<#1>{{\!\left<#1\right>}}

\newtheorem{thm}{Theorem}[section]
\newtheorem{lem}[thm]{Lemma}
\newtheorem{prop}[thm]{Proposition}
\newtheorem{cor}[thm]{Corollary}
\newtheorem{df}[thm]{Definition}
\newtheorem{df-lem}[thm]{Definition-Lemma}

\def\ss{\section}
\def\sss{\subsection}
\def\ssss{\subsubsection}

\def\P{{\mathbb P}}

\def\O{{\mathcal O}}

\def\I{{\mathcal{I}}}

\def\tm{{\times}}
\def\rmk{\medskip\noindent \textbf{Remark: }}
\def\eg{\medskip\noindent \textbf{Example: }}
\def\nid{\noindent}

\subjclass[2000]{14N20, 14C05}
\title{Wonderful compactification of an arrangement of subvarieties}
\author{Li Li}

\begin{document}
\maketitle

\ss {Introduction} The purpose of this paper is to define the
so-called wonderful compactification of an arrangement of
subvarieties, to prove its expected properties, to
give a construction by a sequence of blow-ups and to discuss the
order in which the blow-ups can be carried out.

\medskip

Fix a nonsingular algebraic variety $Y$ over an algebraically closed
field (of arbitrary characteristic). An {\em arrangement} of
subvarieties $\mathcal{S}$ is a finite collection of nonsingular
subvarieties such that all nonempty scheme-theoretic intersections of
subvarieties in $\mathcal{S}$ are again in $\mathcal{S}$, or
equivalently, such that any two subvarieties intersect cleanly and the
intersection is either empty or a subvariety in this collection
(see Definition \ref{def arrangement}).

Let $\mathcal{S}$ be an arrangement of subvarieties of $Y$. A subset
$\mathcal{G}\subseteq \mathcal{S}$ is called a {\em building set
of $\mathcal{S}$} if $\forall S\in \mathcal{S}\setminus
\mathcal{G}$, the minimal elements in $\{G\in\mathcal{G}:G\supseteq S\}$ intersect transversally and the intersection is $S$. A set of subvarieties $\mathcal{G}$ is called a {\em building set}\, if all the possible intersections of subvarieties in $\mathcal{G}$ form an arrangement
$\mathcal{S}$ (called the {\em induced arrangement} of $\mathcal{G}$) and $\mathcal{G}$ is a building set of $\mathcal{S}$ (see Definition \ref{def building set subvariety}).

For any building set $\mathcal{G}$, the wonderful compactification
of $\mathcal{G}$ is defined as follows.

\begin{df}\label{def wonderful compactification} Let $\mathcal{G}$ be a nonempty building
set and $Y^\circ=Y\setminus\cup_{G\in\mathcal{G}}G$. The closure of
the image of the natural locally closed embedding
$$Y^\circ\hookrightarrow \prod_{G\in\mathcal{G}}Bl_GY,$$
is called the {\em wonderful compactification} of the arrangement $\mathcal{G}$ and is denoted
by $Y_\mathcal{G}$.
\end{df}

The following description of $Y_{\mathcal{G}}$ is the main theorem and is proved at the end of \S2.3. A {\em $\mathcal{G}$-nest} is a subset of the building set $\mathcal{G}$ satisfying some inductive condition (see Definition \ref{nest of arrangement}).

\begin{thm}\label{main wonder}
  Let $Y$ be a nonsingular variety and $\mathcal{G}$ be a nonempty building
set of subvarieties of $Y$. Then the wonderful compactification $Y_\mathcal{G}$ is a nonsingular variety. Moreover, for each $G\in\mathcal{G}$
  there is a nonsingular divisor $D_G\subset Y_\mathcal{G}$, such that:
\begin{enumerate}
\item[(i)] The union of these divisors is  $Y_\mathcal{G}\setminus
Y^\circ$.
\item[(ii)] Any set of these divisors meet transversally. An
  intersection of divisors $D_{T_1}\cap\cdots\cap D_{T_r}$ is nonempty
  exactly when $\{T_1,\cdots,T_r\}$ form a $\mathcal{G}$-nest.
\end{enumerate}
\end{thm}

This theorem is proved by a construction of $Y_\mathcal{G}$ through an
explicit sequence of blow-ups of $Y$ along nonsingular centers
(Definition \ref{construction}, Theorem \ref{main prop}).

\medskip
Here are some examples of wonderful compactifications of an arrangement (see \S4 for details).

\begin{enumerate}
\item De Concini-Procesi's wonderful model of subspace arrangements (\S\ref{eg:DP}).
In this case, $Y$ is a vector space, $\mathcal{S}$ is a finite set
of proper subspaces of $Y$ and $\mathcal{G}$ is a building set with
respect to $\mathcal{S}$.

\item
Suppose $X$ is a nonsingular algebraic variety, $n$ is a positive integer and $Y$ is the Cartesian product $X^n$. A diagonal of $X^n$ is $$\Delta_I=\{(p_1,\dots,p_n)\in X^n| p_i=p_j, \forall\, i,j\in I\}$$ for $I\subseteq [n]$, $|I|\ge 2$. A polydiagonal is an intersections of diagonals
$$\Delta_{I_1}\cap\cdots\cap\Delta_{I_k}$$ for $I_i\subseteq [n], |I_i|\ge 2$ ($1\le i\le k$).

\begin{enumerate}
\item Fulton-MacPherson configuration space $X[n]$ (\S\ref{eg:FM}). This is the
wonderful compactification $Y_\mathcal{G}$ where $\mathcal{G}$ is
the set of all diagonals in $Y$ and the induced arrangement
$\mathcal{S}$ is the set of all polydiagonals. It is a special
example of Kuperberg-Thurston's compactification $X^\Gamma$ when
$\Gamma$ is the complete graph with $n$ vertices.

\item Ulyanov's polydiagonal compactification $X\pd<n>$ (\S\ref{eg:U}). It is the
wonderful compactification $Y_\mathcal{G}$ where
$\mathcal{S}=\mathcal{G}$ are the set of all polydiagonals.

\item Kuperberg-Thurston's compactification $X^\Gamma$ where $\Gamma$
is a connected graph with $n$ labeled vertices (\S\ref{eg:KT}).
$X^\Gamma$ is the wonderful compactification $Y_\mathcal{G}$ where
$\mathcal{G}$ is the set of diagonals in $Y$ corresponding to
vertex-2-connected subgraphs of $\Gamma$, and $\mathcal{S}$ is the set
of polydiagonals generated by intersections of diagonals in
$\mathcal{G}$.

\end{enumerate}

\item Moduli space of rational curves with $n$ marked points
$\overline{M}_{0,n}$ (\S\ref{eg:Moduli}). It is the wonderful
compactification $Y_\mathcal{G}$ where $Y=(\mathbb{P}^1)^{n-3}$ and
$\mathcal{G}$ is set of all diagonals and augmented diagonals $\Delta_{I,a}$ defined as
$$\Delta_{I,a}:=\{(p_4,\cdots,p_{n})\in (\mathbb{P}^1)^{n-3}|\, p_i=a, \forall i\in I\}$$ for
$I\subseteq \{4,\dots,n\}$, $|I|\ge 2$ and $a\in \{0, 1, \infty\}$.

The moduli space $\overline{M}_{0,n}$ is also the wonderful compactification $Y_\mathcal{G}$ where
$Y=\mathbb{P}^{n-3}$ and $\mathcal{G}$ is the set of all projective
subspaces of $\mathbb{P}^{n-3}$ spanned by any subset of fixed $n-1$
generic points \cite{Kapranov}.

\item Hu's compactification of open varieties (\S\ref{eg:H}). It is a wonderful
compactification of $(Y,\mathcal{S},\mathcal{G})$ where $Y$ is a
nonsingular algebraic variety, $\mathcal{S}=\mathcal{G}$ is an
arrangement of subvarieties of $Y$ \cite{Hu}.

\end{enumerate}

During the study of the sequence of blow-ups, a natural question
arises: in which order can we carry out the blow-ups to obtain the wonderful
compactification? For example, the original construction of
Fulton-MacPherson configuration space $X[n]$, Keel's construction of
$\overline{M}_{0,n}$, and Kapranov's construction of
$\overline{M}_{0,n}$, none of them are obtained by blowing up along the
centers with increasing dimensions. If we change the order of
blow-ups, do we still get the same variety?

We answer this question with the following theorem which is proved in \S3. The notation $\widetilde{G}$ stands for the so-called {\em dominant tranform} of $G$ (see Definition \ref{sim}) which is similar but slightly different to the strict tranform: for a subvariety $G$ contained in the center of a blow-up, the strict transform of $G$ is empty but the dominant transform $\widetilde{G}$ is the preimage of $G$. 

\begin{thm}\label{order}   Let $Y$ be a nonsingular variety and $\mathcal{G}=\{G_1,\cdots,G_N\}$ be a nonempty building
set of subvarieties of $Y$. Let $\I_i$ be the ideal sheaf of $G_i\in\mathcal{G}$. 
\begin{itemize}
\item[(i)] The wonderful compactification $Y_{\mathcal{G}}$ is isomorphic
to the blow-up of \,$Y$ along the ideal sheaf \,$\I_1\I_2\cdots \I_N$.

\item[(ii)] If we arrange $\mathcal{G}=\{G_1,\dots, G_N\}$ in such an order that
$$\textrm{\emph{$(*)$\, the first
$i$ terms $G_1,\dots, G_i$ form a building set for any $1\le i\le N$.}}$$
then
$$Y_\mathcal{G}=Bl_{\widetilde{G}_N}\cdots Bl_{\widetilde{G}_2}Bl_{G_1}Y,$$ where each blow-up is
along a nonsingular subvariety.
\end{itemize}
\end{thm}

\eg By Keel's construction \cite{Keel1} and the above theorem,
$\overline{M}_{0,n}$ is isomorphic to the wonderful compactification 
$Y_\mathcal{G}$ where $Y$ is $(\mathbb{P}^1)^{n-3}$ and
$\mathcal{G}$ is set of all diagonals and augmented
diagonals. In other words, we can blow up along the centers in any order satisfying
(*) (e.g. of increasing dimension). As a consequence, we have

\begin{cor}\label{corM0n}
 Let $\psi: \P^1[n]\to (\P^1)^3$ be the composition of the natural morphism $\P^1[n]\to (\P^1)^n$ and
 let $\pi_{123}:( \P^1)^n\to (\P^1)^3$ be the projection to the first three components. Then $\overline{M}_{0,n}$ is isomorphic to the fiber of $\psi$ over the point
 $(0,1,\infty)\in(\P^1)^3$. Equivalently, $\overline{M}_{0,n}$ is isomorphic to  the fiber over any point $(p_1,p_2,p_3)$ where
 $p_1,p_2,p_3$ are three distinct points in $\P^1$.
\end{cor}

Similarly, Kapranov's construction does not delicately depend on
the order of the blow-ups, for example we can blow up along the centers in any order of
increasing dimension.

\medskip

This article is built on the following previous works:
Fulton-MacPherson \cite{FM}, De Concini-Procesi \cite{DP},
MacPherson-Procesi \cite{MP}, Ulyanov \cite{Ulyanov}, Hu \cite{Hu}.

The inspiring paper by De Concini and Procesi \cite{DP} gives a thorough discussion of an
arrangement of linear subspaces of a vector space. Given a vector space $Y$ and an arrangement of
subspaces $\mathcal{S}$, De Concini and Procesi give a condition for a subset
$\mathcal{G}\subseteq\mathcal{S}$ such that there exists a so called \emph{wonderful model}
$Y_\mathcal{G}$ of the arrangement, in which the elements in $\mathcal{G}$ are replaced by simple
normal crossing divisors. De Concini and Procesi call $\mathcal{G}$ a \emph{building set}.
Their paper also gives a criterion of whether the intersection of a collection of such divisors is
nonempty by introducing the notion of a \emph{nest}.

Later, this idea has been generalized to nonsingular varieties over $\mathbb{C}$ with conical
stratifications by MacPherson and Procesi. They consider conical stratifications in place of subspace arrangments in \cite{DP}. The notion of \emph{building set} and \emph{nest} is
generalized in this setting. The idea of the construction of wonderful compactifications of
arrangement of subvarieties in our paper is largely inspired by the beautiful paper \cite{MP}.

In our paper, we give definitions of arrangements of
subvarieties, building sets and nests. The wonderful
compactifications are shown to have properties analogous to the ones
in \cite{DP} or \cite{MP}.

The paper is organized as follows. In section 2 we give the
construction of the wonderful compactification $Y_\mathcal{G}$. In
\S2.1 we give the definition of arrangements, building sets and
nests. \S2.2 is the description of how the arrangements, building
sets and nests vary under one blow-up. \S2.3 gives the construction
of $Y_\mathcal{G}$. In section 3 we discuss that in which order could the blow-ups be carried out to obtain $Y_\mathcal{G}$.
Section 4 gives some examples of wonderful compactifications. In \S5.1 we discuss clean intersections and transversal
intersections. In \S5.2 we give the proofs of previous statements. In \S5.3 we discuss how different choices of blow-ups change the codimension of the centers. Finally in \S5.4, we give the statements for a general (non-simple) arrangement (proofs omitted).

\noindent\textbf{Acknowledgements.} In many ways the author greatly
indebted to Mark de Cataldo, his Ph.D. advisor. He is very grateful to
William Fulton for valuable comments. He would also thank Blaine
Lawson, Dror Varolin, Jun-Muk Hwang, especially Herwig Hauser, for
their many useful comments and encouragement. He thanks Jonah Sinick for carefully proof-reading the paper. He thanks the referee for many constructive suggestions to improve the presentation.

\ss {Arrangements of subvarieties and the wonderful
compactifications}

By a variety we shall mean a reduced and irreducible
algebraic scheme defined over a fixed algebraically closed field (of
arbitrary characteristic). A subvariety of a variety is a closed
subscheme which is a variety. By a point of a variety we shall mean a
closed point of that variety. By the intersection of subvarieties
$Z_1, \dots, Z_k$ we shall mean the set-theoretic intersection (denoted by
$Z_1\cap\dots\cap Z_k$). We denote the ideal sheaf of a subvariety
$V$ of a variety $Y$ by $\I_V$.

In this section, we will discuss arrangements, building sets and nests, based on which we define the wonderful compactifications of an arrangement.
The idea is inspired by \cite{DP}\cite{MP}.

\sss{Arrangement, building set, nest} The following definition of arrangement is adapted from \cite{Hu}.  For a brief review of the
definitions of clean intersection and of transversal intersection,
please see Appendix \S\ref{intersections}.
\begin{df}\label{def arrangement}
A simple {\em arrangement} of subvarieties of a nonsingular variety\, $Y$ is a finite set $\mathcal{S}=\{S_i\}$ of nonsingular closed
subvarieties $S_i$ properly contained in $Y$ satisfying the following conditions:
\begin{itemize}
\item[(i)] $S_i$ and $S_j$ intersect cleanly (i.e. their intersection is nonsingular and the tangent
bundles satisfy $T(S_i\cap S_j)=T(S_i)|_{(S_i\cap S_j)}\cap
T(S_j)|_{(S_i\cap S_j)}$),
\item[(ii)] $S_i\cap S_j$ either equals to some $S_k$ or is empty.
\end{itemize}
\end{df}

The above definition is equivalent to say that $\mathcal{S}$ is an arrangement if and only if it is closed under scheme-theoretic intersections, cf. Lemma \ref{ideal_clean}.

 Although we will discuss only the simple arrangement for
simplicity, most statements still hold, with minor revision, for
general arrangements, i.e., instead of the above condition (2), we
allow $S_i\cap S_j$ to be a disjoint union of some $S_k$'s. (See Appendix \S\ref{general}.)

For a simple arrangement, the condition of transversality can be
checked at one point (instead of at every point) of the intersection
(Lemma \ref{appendix lemma}).

\begin{df}\label{def building set subvariety}
 Let $\mathcal{S}$ be an arrangement of subvarieties of $\, Y$. A subset $\mathcal{G}\subseteq \mathcal{S}$ is
called a {\em building set of $\mathcal{S}$} if $\forall S\in \mathcal{S}$, the minimal elements in $\{G\in\mathcal{G}: G\supseteq S\}$ intersect
transversally and their intersection is $S$ (by our definition of transversality \S\ref{intersections}, the condition is satisfied if $S\in
\mathcal{G}$). In this case, these minimal elements are called the {\em $\mathcal{G}$-factors} of $S$.

A finite set $\mathcal{G}$ of nonsingular subvarieties of $Y$ is
called a {\em building set} if the set of all possible
intersections of collections of subvarieties from $\mathcal{G}$
forms an arrangement $\mathcal{S}$, and that $\mathcal{G}$ is a
building set of $\mathcal{S}$. In this situation,
$\mathcal{S}$ is called the {\em arrangement induced by
$\mathcal{G}$}.
\end{df}

\eg Let $X$ be a nonsingular variety of positive dimension and $Y$ be the Cartesian product $X^3$. 

\begin{enumerate}

\item The set $\mathcal{G}=\{\Delta_{12},\Delta_{13},\Delta_{23},\Delta_{123}\}$ is a building set whose induced arrangement is
$\mathcal{G}$ itself.
\item The set $\mathcal{G}=\{\Delta_{12},\Delta_{13}\}$ is a building set whose induced arrangement is $\{\Delta_{12},\Delta_{13},\Delta_{123}\}$. On the other hand, $\mathcal{G}$ is not a building set of the arrangement $\{\Delta_{12},\Delta_{13},\Delta_{23},\Delta_{123}\}$.

\item The set $\mathcal{G}=\{\Delta_{12},\Delta_{13},\Delta_{23}\}$ is not a building set, since the set of all possible intersections from
$\mathcal{G}$ is $\{\Delta_{12},\Delta_{13},\Delta_{23},\Delta_{123}\}$, but $\Delta_{123}$ is not a transversal intersection of
$\Delta_{12},\Delta_{13}$ and $\Delta_{23}$.
\end{enumerate}

\rmk The building set $\mathcal{G}$ defined here is related to the one defined in \cite{DP} as follows. For any point $y\in Y$, define
$\mathcal{S}^*_y=\{T_{S,y}^\perp\}_{S\in\mathcal{S}}$ and $\mathcal{G}^*_y=\{T_{S,y}^\perp\}_{S\in\mathcal{G}}$. We claim that the set $\mathcal{G}$ is a building set if and only if $\mathcal{G}^*_y$ is a building set for all $y\in Y$ in the sense of DeConcini and Procesi.

Indeed, $\mathcal{S}$ being an arrangement is equivalent to the
condition that for any $y\in Y$, $\mathcal{S}^*_y$ is a finite set
of nonzero linear subspaces of $T_y^*$ which is closed under sum,
and such that each element of $\mathcal{S}^*_y$ is equal to
$T_{S,y}^\perp$ for a unique $S\in\mathcal{S}$. The subset $\mathcal{G}\subseteq \mathcal{S}$ being a building set is equivalent to the
condition that $\forall S\in \mathcal{S}$, $\forall y\in S$, suppose
$T_1^\perp,\dots, T_k^\perp$ are all the maximal elements of
$\mathcal{G}^*_y$ contained in $T_{S,y}^\perp$, then they form a
direct sum and
$$T_1^\perp\oplus T_2^\perp\oplus\cdots\oplus T_k^\perp=T_{S,y}^\perp,$$
which is exactly the definition of building set in \cite{DP} \S2.3
Theorem (2).

\begin{df}(cf. \cite{MP} \S4)\label{nest of arrangement}
A subset $\mathcal{T}\subseteq\mathcal{G}$ is called
{\em $\mathcal{G}$-nested} (or a {\em $\mathcal{G}$-nest}) if
it satisfies one of the following equivalent conditions:
\begin{enumerate}
\item[(i)] There is a flag of elements in $\mathcal{S}$: $S_1\subseteq S_2\subseteq\dots\subseteq S_\ell$, such that
 $$\mathcal{T}=\bigcup_{i=1}^\ell\{A: \hbox{ $A$ is a $\mathcal{G}$-factor of $S_i$}\}.$$
 (We say $\mathcal{T}$ is {\em induced} by the flag $S_1\subseteq S_2\subseteq\dots\subseteq S_\ell$.)
\item[(ii)] Let $A_1,\dots,A_k$ be the minimal elements of $\mathcal{T}$, then they are all the $\mathcal{G}$-factors of
certain element in $\mathcal{S}$. For any $1\le i\le k$, the set $\{A\in \mathcal{T}: A\supsetneq A_i\}$ is also $\mathcal{G}$-nested
defined by induction.
\end{enumerate}
\end{df}

\eg Let $X$ be a nonsingular variety of positive dimension and $Y$ be the Cartesian product $X^4$. Take the building set $\mathcal{G}$ to be the set of all diagonals in $X^4$.

\begin{enumerate}
 
\item  The set 
$\mathcal{T}=\{\Delta_{12}, \Delta_{123}\}$ is a $\mathcal{G}$-nest, since it can be induced by the flag $\Delta_{123}\subseteq\Delta_{12}$.

\item The set 
$\mathcal{T}=\{\Delta_{12}, \Delta_{34}, \Delta_{1234}\}$ is a $\mathcal{G}$-nest, since it can be induced by the flag $\Delta_{1234}\subseteq(\Delta_{12}\cap \Delta_{34})$.

\item The set $\mathcal{T}=\{\Delta_{12},\Delta_{13}\}$ is not a $\mathcal{G}$-nest. Indeed, the intersection of the minimal elements in $\mathcal{T}$ is $\Delta_{123}$, which has only one $\mathcal{G}$-factor: $\Delta_{123}$ itself. By condition (ii) of the definition, $\mathcal{T}$ is not a $\mathcal{G}$-nest.
\end{enumerate}

Note that the intersection of elements in a
$\mathcal{G}$-nest $\mathcal{T}$ is nonempty by (ii). Now we explain why the two conditions (i) and (ii) are equivalent.
Given a set $\mathcal{T}$ satisfying (ii), we can construct a flag as follows: define $S_1=A_1\cap\cdots\cap A_k$ which is the intersection of all
subvarieties in $\mathcal{T}$. Let $S_2$ be the intersection of the subvarieties in $\mathcal{T}$ which are not minimal elements in $\mathcal{T}$
that contain $S_1$. Then inductively let $S_{j+1}$ be the intersection of those which are not minimal elements in $\mathcal{T}$ that contains $S_j$. It
is easy to show that $\mathcal{T}$ is induced by the flag $S_1\subseteq S_2\subseteq\cdots$,  therefore (ii)$\Rightarrow$(i). On the other hand, let
$S_{1j}=A_j$ ($1\le j\le k$) be the $\mathcal{G}$-factors of $S_1$. Note that for any $1\le i\le \ell$, a $\mathcal{G}$-factor of $S_i$ must contain
exactly one element of $A_1,\dots,A_k$, otherwise the $A_i$'s will not intersect transversally. Let $S_{ij}$ be the $\mathcal{G}$-factor of $S_i$
that contains $A_j$. (Define $S_{ij}=Y$ if there is no such a $\mathcal{G}$-factor.) Then for each $1\le j\le k$, there is a flag of elements
$S_{2j}\subseteq S_{3j}\subseteq\cdots\subseteq S_{\ell j}$, which induces the $\mathcal{G}$-nest $\{A\in \mathcal{T}: A\supsetneq A_i\}$. This shows
(i)$\Rightarrow$(ii).

\medskip
We state some basic properties about arrangements and building sets.

\begin{lem}\label{fact 2}Let $Y$ be a nonsingular variety and $\mathcal{G}$ be a building set with the induced arrangement $\mathcal{S}$. Suppose $S\in\mathcal{S}$ and $G_1,\dots,G_k$ are all the
$\mathcal{G}$-factors of $S$ (Definition \ref{def building set subvariety}). Then

\begin{itemize}
\item[(i)] For any $1\le m\le k$, the subvarieties $G_1,\dots, G_m$ are all the $\mathcal{G}$-factors of
the subvariety $G_1\cap\cdots\cap G_m$.

\item[(ii)] Suppose $F\in\mathcal{G}$ is minimal such that $F\cap S\neq\emptyset$, $F\subseteq
G_1,\dots, G_m$ and $F\nsubseteq G_{m+1},\dots, G_k$. Then $F, G_{m+1},\dots, G_k$ are all the
$\mathcal{G}$-factors of the subvariety $F\cap S$.
\end{itemize}
\end{lem}

\begin{proof}
  See Appendix \S5.2.
\end{proof}

Here is an immediate consequence of Lemma \ref{fact 2}:

\begin{lem}\label{minimals}
If $G_1,\dots, G_k\in \mathcal{G}$ are all minimal and their intersection $S$ is nonempty, then
$G_1,\cdots, G_k$ are all the $\mathcal{G}$-factors of $S$.
\end{lem}

Next, we introduce the notion of the $F$-factorization, which turns
out to be a convenient terminology for the proof of the construction of
wonderful compactifications.

\begin{df-lem}\label{fact 1}

 Suppose $F\in\mathcal{G}$ is minimal. Then
\begin{itemize}
\item[(i)] Any $G\in\mathcal{G}$ either contains $F$ or intersects transversally with $F$.

\item[(ii)] Every $S\in\mathcal{S}$ satisfying $S\cap F\neq\emptyset$ can be uniquely expressed as
 $A\cap B$ where $A, B\in\mathcal{S}\cup\{Y\}$ satisfy $A\supseteq F$ and $B\pitchfork F$ (hence $A\pitchfork B$).
 We call this expression $S=A\cap B$ the {\em $F$-factorization} of $S$.
\item[(iii)] Suppose the $\mathcal{G}$-factors of $S$ are $G_1, \dots, G_k$ where $G_1,\dots, G_m$
contain $F$ ($0\le m\le k$, the case $m=0$ is understood to mean
that no $\mathcal{G}$-factors of $S$ contain $F$) and let the
$F$-factorization of $S$ be $A\cap B$.

Then $G_1,\dots, G_m$ are all the $\mathcal{G}$-factors of $A$ and
$G_{m+1},\dots, G_k$ are all the $\mathcal{G}$-factors of $B$, so
$A=\cap_{i=1}^mG_i$ and $B=\cap_{i=m+1}^kG_i$. (Here we assume $A=Y$
if $m=0$, assume $B=Y$ if $m=k$.)

\item[(iv)] Suppose $S'\in\mathcal{S}$ such that $S'\cap S\cap F\neq\emptyset$.  Let $S'=A'\cap
B'$ be the $F$-factorization of $S'$.

Then $F\pitchfork (B\cap B')$, therefore the $F$-factorization of $S\cap S'$ is $(A\cap A')\cap
(B\cap B')$.
\end{itemize}
\end{df-lem}
\begin{proof}
  See Appendix \S5.2.
\end{proof}

\sss{Change of an arrangement after a blow-up}

Before considering a sequence of blow-ups, we first consider a single
blow-up. Let $Y$ be a nonsingular variety and $\mathcal{G}$ be a building set with the induced arrangement $\mathcal{S}$. In Proposition \ref{arrangment in blow-up} we show that, if $F\in\mathcal{G}$ is minimal, then
there exists a natural arrangement $\widetilde{\mathcal{S}}$ in
$Bl_FY$ induced from $\mathcal{S}$ and a natural building set
$\widetilde{\mathcal{G}}$  induced from $\mathcal{G}$.

\begin{df}\label{sim} Let $Z$ be a nonsingular subvariety of a nonsingular variety $Y$ and $\pi: Bl_ZY\to Y$ be the blow-up of $Y$ along $Z$.
 
For any irreducible subvariety $V$ of $Y$, we define the {\em dominant transform} of $V$, denoted by $\widetilde{V}$ or $V^\sim$, to be the strict transform of $V$ if
$V\nsubseteq G$ and to be the scheme-theoretic inverse $\pi^{-1}(V)$ if
$V\subseteq G$.

For a sequence of blow-ups, we still denote the
iterated dominant transform $(\cdots((V^\sim)^\sim)\cdots)^\sim$ by $\widetilde{V}$ or $V^\sim$.
\end{df}
\rmk The reason that we introduce the notion of dominant transform is because that the strict transform does not behave as expected: the strict transform of a subvariety contained in the center of a blow-up is empty, which is not what we need. 

\begin{prop}\label{arrangment in blow-up} Let $Y$ be a nonsingular variety and $\mathcal{G}$ be a building set with the induced arrangement $\mathcal{S}$.  Let $F$ be a minimal
element in $\mathcal{G}$ and let $\pi: Bl_FY\to Y$ be the blow-up of
\,$Y$ along $F$. Denote the exceptional divisor by $E$.
\begin{itemize}
\item[(i)] The collection $\widetilde{\mathcal{S}}$ of subvarieties in $Bl_GY$ defined as
  $$\widetilde{\mathcal{S}}:=\{\widetilde{S}\}_{S\in\mathcal{S}}\cup\{\widetilde{S}\cap E\}_
  {\emptyset\subsetneq S\cap F\subsetneq S}$$ is a (simple) arrangement of subvarieties in $Bl_GY$.

\item[(ii)] $\widetilde{\mathcal{G}}:=\{\widetilde{G}\}_{G\in\mathcal{G}}$ is a building set
of $\widetilde{\mathcal{S}}$.

\item[(iii)] Given a subset $\mathcal{T}$ of $\mathcal{G}$, we define $\widetilde{\mathcal{T}}:=
  \{\widetilde{A}\}_{A\in\mathcal{T}}$. Then $\mathcal{T}$ is a $\mathcal{G}$-nest if and only if
  $\widetilde{\mathcal{T}}$ is a $\widetilde{\mathcal{G}}$-nest.
\end{itemize}
\end{prop}

The proof is in Appendix \S5.2. The main ingredient is the following lemma.

\begin{lem}\label{basic} Assume the same notation as in Proposition \ref{arrangment in blow-up}.  Assume $A$, $A_1$, $A_2$, $B$, $B_1$, $B_2$, $G$  are nonsingular subvarieties of
$Y$. 
\begin{itemize}
\item[(i)] Suppose $A\supsetneq F$. Then $\widetilde{A}\cap E$
intersect transversally (hence cleanly).
\item[(ii)] Suppose $A_1\nsubseteq A_2$,
$A_2\nsubseteq A_1$ and $A_1\cap A_2=F$ and the intersection is
clean. Then $\widetilde{A}_1\cap\widetilde{A}_2=\emptyset$.
\item[(iii)]  Suppose $A_1$ and $A_2$ intersect cleanly and $F\subsetneq A_1\cap A_2$. Then
$\widetilde{A}_1\cap\widetilde{A}_2={(A_1\cap
  A_2)}^\sim$. Moreover, $\widetilde{A}_1$ and $\widetilde{A}_2$ intersect cleanly.
\item[(iv)] Suppose $B_1$ and $B_2$ intersect cleanly, and $G$ is transversal to $B_1$, $B_2$ and $B_1\cap B_2$. Then $\widetilde{B}_1\cap\widetilde{B}_2={(B_1\cap
  B_2)}^\sim$. Moreover, $\widetilde{B}_1$ and $\widetilde{B}_2$ intersect cleanly.
\item[(v)] Suppose $A\pitchfork B$, $F\subseteq A$ and $F\pitchfork B$. Then $\widetilde{A}\cap\widetilde{B}=
  (A\cap B)^\sim$. Moreover, $\widetilde{A}\pitchfork\widetilde{B}$ and $(E\cap \widetilde{A})\pitchfork \widetilde{B}.$
\item[(vi)] Assume $F\subseteq A$, $F\pitchfork B_1\pitchfork B_2$, $G\subseteq F\cap B_1$ and
$G\pitchfork B_2$. Then $\widetilde{G}\cap \widetilde{A}\cap (B_1\cap
B_2)^\sim=\widetilde{G}\cap \widetilde{A}\cap \widetilde{B}_2$ where the latter
is a transversal intersection.
\end{itemize}
\end{lem}

\begin{proof}
  See Appendix \S5.2.
\end{proof}

\medskip
\sss{A sequence of blow-ups and the construction of wonderful
compactifications} Now we study a sequence of blow-ups, give
different descriptions of a wonderful compactification, and study
the relations of the arrangements occurred in the sequence of
blow-ups.

Given $k$ morphisms between algebraic varieties with the same domain
$f_i: X\to Y_i$, we adopt the notation $(f_1, f_2, \cdots,
f_k): X\to Y_1\times\cdots\times Y_k$ to be the composition of the
diagonal morphism $X\to X\times\cdots\times X$ with the morphism
$f_1\times\cdots\times f_k$.
\begin{lem}\label{blowup embedding}
Let $V$ and $W$ be two nonsingular
  subvarieties of a nonsingular variety $Y$ such that either $V$ and $W$ intersect transversally or one of $V$ and $W$ contains the other. Let $f: Y_1\to Y$
  (resp. $g: Y_2\to Y$) be the blow-up of\, $Y$ along $W$ (resp. $V$). Let
  $g': Y_3\to Y_1$ be the blow-up of $Y_1$ along the dominant transform $\widetilde{V}$. Then there exists a morphism $f':
  Y_3\to Y_2$ such that the following diagram commutes:
$$\xymatrix{Y_3\ar[r]^{f'}\ar[d]^{g'} & Y_2\ar[d]^{g}\\ Y_1\ar[r]^{f} & Y}$$
Moreover, $(g',f'): Y_3\to Y_1\tm Y_2$ is a closed embedding.

\end{lem}

\begin{proof}
  Because of the universal property of blowing up (\cite{Ha} Proposition 7.14), to show the existence
  of $f'$, we need only to show that $(fg')^{-1}\I_V\cdot \O_{Y_3}$ is an invertible
  sheaf of ideals on $Y_3$. This is true because by our choice of $V$ and $W$, the sheaf $f^{-1}\I_V\cdot \O_{Y_1}$ is either $\I_{\widetilde{V}}$ or $\I_{\widetilde{V}}\I_E$, where $E$ is the exceptional divisor of the blow-up $f: Y_1\to
  Y$. Hence the ideal sheaf
 $$
  (fg')^{-1}\I_V\cdot \O_{Y_3}=g'^{-1}(f^{-1}\I_V\cdot \O_{Y_1})\cdot
  \O_{Y_3}$$ is either $(g'^{-1}\I_{\widetilde{V}})\cdot \O_{Y_3}$
  or
  $g'^{-1}(\I_{\widetilde{V}}\cdot \I_E)\cdot \O_{Y_3}$, both of which are invertible by the construction of $g'$, therefore the ideal sheaf $(fg')^{-1}\I_V\cdot \O_{Y_3}$ is invertible.

  The fact that $(g',f')$ is a closed embedding can be checked using local parameters.
\end{proof}

\begin{lem}\label{embedding induction}
  Suppose $X_1, X_2, X_3, Y_1, Y_2, Y_3$ are nonsingular varieties such that the following diagram
  commutes,
  $$\xymatrix{X_1\ar[r]^{f_1}\ar[d]^{g_1} &X_2\ar[r]^{f_2}\ar[d]^{g_2} &X_3\ar[d]^{g_3}\\
     Y_1\ar[r]^{h_1} &Y_2\ar[r]^{h_2} &Y_3}$$
  If $(g_1,f_1): X_1\to Y_1\tm X_2$ and $(g_2,f_2): X_2\to Y_2\tm X_3$ are closed embeddings,
  then $(g_1,f_2f_1): X_1\to Y_1\tm X_3$ is also a closed embedding.

  As a consequence, if we have the following commutative diagram
  $$\xymatrix{X_1\ar[r]^{f_1}\ar[d]^{g_1} &X_2\ar[r]^{f_2}\ar[d]^{g_2} &\cdots\ar[r]^{f_{k-1}} &X_k\ar[d]^{g_k}\\
     Y_1\ar[r]^{h_1} &Y_2\ar[r]^{h_2} &\cdots\ar[r]^{h_{k-1}} &Y_k}$$
  and $(g_i,f_i): X_i\to Y_i\tm X_{i+1}$ are closed embeddings for all $1\le i\le k-1$, then
  $(g_1,f_{k-1}\cdots f_1): X_1\to Y_1\tm X_k$ is also a closed embedding.
\end{lem}
\begin{proof}
  The composition of two closed embeddings is still a closed embedding, so
  $$\phi:=(g_1, g_2f_1, f_2f_1): X_1\to Y_1\tm Y_2\tm X_3$$
  is a closed embedding, whose image $\phi(X_1)$ is a closed subvariety of $Y_1\tm Y_2\tm
  X_3$ which is isomorphic to $X_1$.  Consider the projection $\pi_{13}: Y_1\tm Y_2\tm  X_3\to Y_1\tm X_3$, and the morphism
  $\Gamma_{h_1}\tm 1_{X_3}: Y_1\tm X_3\to Y_1\tm Y_2\tm X_3$.
$$\xymatrix{X_1\ar[r]^-{\phi}\ar[rd]_{(g_1,f_2f_1)}&Y_1\times Y_2\times X_3\ar@<-.5ex>[d]_{\pi_{13}}\\ & Y_1\times X_3\ar@<-.5ex>[u]_{\Gamma_{h_1}\times 1_{X_3}}}$$
  Notice that $\pi_{13}\circ (\Gamma_{h_1}\tm
  1_{X_3})$ is the identity automorphism of $Y_1\tm X_3$, and
  $(\Gamma_{h_1}\tm 1_{X_3})\circ\pi_{13}|_{\phi(X_1)}$ is the identity automorphism of $\phi(X_1)$.
  It follows that $(g_1,f_2f_1): X_1\to Y_1\tm X_3$ is a closed embedding.
\end{proof}

\begin{df}[Inductive construction of $Y_\mathcal{G}$]\label{construction} Let $Y$ be a nonsingular variety, $\mathcal{S}$ be an arrangement of
subvarieties and $\mathcal{G}$ be a building set of $\mathcal{S}$.
Suppose $\mathcal{G}=\{G_1,\dots, G_N\}$ is indexed in an order
compatible with inclusion relations, i.e. $i\le j $ if $
G_i\subseteq G_j$. We define $(Y_k, \mathcal{S}^{(k)},
\mathcal{G}^{(k)})$ inductively with respect to $k$:

(i) For $k=0$, define $Y_0=Y$, $\mathcal{S}^{(0)}=\mathcal{S}$,
$\mathcal{G}^{(0)}=\mathcal{G}=\{G_1,\dots, G_N\}$, $G^{(0)}_i=G_i$ for $1\le i\le N$.

(ii) Assume that $(Y_{k-1},\mathcal{S}^{(k-1)},
\mathcal{G}^{(k-1)})$ is constructed. 
\begin{itemize}
 \item 
Define $Y_{k}$ to be the
blow-up of $Y_{k-1}$ along the nonsingular subvariety $G^{(k-1)}_k$.
\item Define $G^{(k)}:=(G^{(k-1)})^\sim$ for $G\in \mathcal{G}$ and define
$$\mathcal{G}^{(k)}:=\{G^{(k)}\}_{G\in\mathcal{G}}.$$
\item Define $\mathcal{S}^{(k)}$ to be the induced arrangement of $\mathcal{G}^{(k)}$.
\end{itemize}

(iii) Continue the inductive construction until $k=N$. We obtain $$(Y_N, \mathcal{S}^{(N)},
\mathcal{G}^{(N)})$$ where all the subvarieties in the building set
$\mathcal{G}^{(N)}$ are divisors.
\end{df}

\rmk In step (ii) we need Proposition \ref{arrangment in blow-up}. Indeed, since $G^{(k-1)}_i$ for $i<k$ are all divisors hence are too large
to be contained in $G^{(k-1)}_k$, $G^{(k-1)}_k$ is minimal in
$\mathcal{G}^{(k-1)}$. Proposition \ref{arrangment in blow-up} then
asserts the existence of a naturally induced arrangement
$\mathcal{S}^{(k)}$, and asserts that
$\mathcal{G}^{(k)}=\{G^{(k)}\}_{G\in\mathcal{G}}$ is a building set
with respect to $\mathcal{S}^{(k)}$. 

\begin{prop}\label{main prop}
  The variety $Y_N$ constructed in Definition \ref{construction} is isomorphic to the wonderful compactification $Y_\mathcal{G}$ defined in Definition \ref{def wonderful compactification}.
\end{prop}

\begin{proof}
  We prove by induction that $Y_k$ is the closure of the inclusion $$Y^\circ\hookrightarrow \prod_{i=1}^k
  Bl_{G_i}Y.$$  The proposition is then the special case $k=N$.

  Let $0\le i\le k-1$. Since $G^{(i)}_{i+1}$ is minimal in $\mathcal{G}^{(i)}$, Definition-Lemma \ref{fact 1} (i)
  asserts that there are only two possible relations between the nonsingular
  subvarieties $G^{(i)}_k$ and $G^{(i)}_{i+1}$ of $Y_i$: either $G^{(i)}_k\supseteq G^{(i)}_{i+1}$
  or $G^{(i)}_k\pitchfork G^{(i)}_{i+1}$.
  Therefore Lemma \ref{blowup embedding} applies. Since $G^{(i+1)}_k=(G^{(i)}_k)^\sim$,  there exists a
  morphism $f'$ such that following diagram commutes.
  $$\xymatrix{Bl_{G^{(i+1)}_k}Y_{i+1}\ar[r]^{f'}\ar[d]^{g'} & Bl_{G^{(i)}_k}Y_{i}\ar[d]^{g}\\ Y_{i+1}\ar[r]^{f} & Y_i}$$
  The morphism $(g', f'): Bl_{G^{(i+1)}_k}Y_{i+1}\to Y_{i+1}\tm Bl_{G^{(i)}_k}Y_{i}$ is a closed
  embedding. Using Lemma \ref{embedding induction} on the following diagram
   $$\xymatrix{Bl_{G^{(k-1)}_k}Y_{k-1}\ar[r]\ar[d] &Bl_{G^{(k-2)}_k}Y_{k-2}\ar[r]\ar[d] &\cdots\ar[r] &Bl_{G^{(0)}_k}Y_{0}\ar[d]\\
     Y_{k-1}\ar[r] &Y_{k-2}\ar[r] &\cdots\ar[r] &Y_0}$$
  and the facts that $Y_k=Bl_{G^{(k-1)}_k}Y_{k-1}$, $G^{(0)}_k=G_k$ and $Y_0=Y$, we conclude that the morphism $$Y_k\to Y_{k-1}\tm Bl_{G_k}Y$$ is a closed
  embedding. Because the composition of closed embeddings is still a closed embedding, the morphism $$Y_k\to
  \prod^k_{i=1}Bl_{G_i}Y$$ is a closed embedding. Then since $Y^{\circ}$ is an open subset of $Y_k$ and since $Y_k$ is irreducible, from the following composition $$Y^\circ\hookrightarrow
   Y_k\hookrightarrow  Y\tm\prod^k_{i=1}Bl_{G_i}Y$$ we see that the closure of $Y^{\circ}$ in $Y\tm\prod^k_{i=1}Bl_{G_i}Y$ is
   $Y_k$.
\end{proof}

\begin{proof}[Proof of Theorem \ref{main wonder}]
Since $Y_\mathcal{G}\cong Y_N$, $Y_\mathcal{G}$ is nonsingular,
$D_G:=G^{(N)}$ are codimension one nonsingular subvarieties of $Y_\mathcal{G}$ and $Y_\mathcal{G}\setminus Y^\circ=\cup D_G$. Therefore
(i) is clear.

For any $T_1, \dots, T_r$ in $\mathcal{G}$ that form a
$\mathcal{G}$-nest, $D_{T_1},\dots,D_{T_r}$ form a
$\mathcal{G}^{(N)}$-nest, hence
$$D_{T_1}\cap\dots\cap D_{T_r}\neq\emptyset$$ by the definition of nest. Conversely, given
 $T_1, \dots, T_r$ in $\mathcal{G}$ such that the above intersection is nonempty, Lemma
 \ref{minimals} implies that
 $D_{T_1},\dots, D_{T_r}$ are all the $\mathcal{G}^{(N)}$-factors of the intersection and therefore intersect transversally.
 Moreover, by the definition of nest,
 $D_{T_1},\dots,D_{T_r}$ form a $\mathcal{G}^{(N)}$-nest. Proposition \ref{arrangment in blow-up} then implies that
 $T_1, \dots, T_r$ form a $\mathcal{G}$-nest. So (ii) is clear.
\end{proof}


\ss{Order of blow-ups} In this section we shall prove Theorem \ref{order}. We shall use this theorem in \S\ref{eg:FM}, \S\ref{eg:KT} and \S\ref{eg:Moduli}. For the proof, we need the following proposition which is stronger than
  Proposition \ref{arrangment in blow-up} (2) in the sense that
  a building set still induces a building set after a blow-up even when the center  of the
  blow-up is not assumed to be minimal.
\begin{prop}\label{G+}
  Suppose $\mathcal{G}=\{G_1,\dots,G_k\}$ is a building set of an arrangement $\mathcal{S}$ in $Y$ and $F\in\mathcal{G}$ is minimal.  Let $\phi: Bl_FY\to Y$ be the blow-up of $Y$ along $F$,  let $\widetilde{\mathcal{G}}$ be
  the induced building set and let $\widetilde{\mathcal{S}}$ be the arrangement induced by $\widetilde{\mathcal{G}}$.
Suppose
  $\mathcal{G}_+=\{G_0,\dots,G_k\}$ is a building set and
  $\mathcal{S}_+$ is the arrangement induced by $\mathcal{G}_+$.

  Then $\widetilde{\mathcal{G}}_+:=\widetilde{\mathcal{G}}\cup\{\widetilde{G}_0\}$ is a building set of the induced arrangement
  $$\widetilde{\mathcal{S}}_+:=\widetilde{\mathcal{S}}\cup \{\widetilde{S}\cap\widetilde{G}_0\}_{S\in\mathcal{S}}.$$
\end{prop}
\begin{proof}
  As the proof of Proposition \ref{arrangment in blow-up}, we need to discuss different types of
intersections of subvarieties. See Appendix
\S5.2.
\end{proof}

\begin{lem}\label{blowup ideal}
  Let $\mathcal{I}_1, \mathcal{I}_2$ be two ideal sheaves on a variety $Y$.  Define $Bl_{\mathcal{I}_2}Bl_{\mathcal{I}_1}Y$ to be the blow-up of \,$Y'=Bl_{\mathcal{I}_1}Y$ along the ideal sheaf $\phi^{-1}\mathcal{I}_2\cdot\mathcal{O}_{Y'}$  where $\phi$ is the blow-up morphism $\phi: Y'\to Y$. Define  $Bl_{\mathcal{I}_1}Bl_{\mathcal{I}_2}Y$ symmetrically.
Then
  $$Bl_{\mathcal{I}_1\mathcal{I}_2}Y\cong Bl_{\mathcal{I}_2}Bl_{\mathcal{I}_1}Y\cong Bl_{\mathcal{I}_1}Bl_{\mathcal{I}_2}Y.$$
 
\end{lem}
\begin{proof}
We show the existence of two natural morphisms
\begin{align*}
&f: Bl_{\mathcal{I}_2}Bl_{\mathcal{I}_1}Y\to Bl_{\mathcal{I}_1\mathcal{I}_2}Y,\\
&g: Bl_{\mathcal{I}_1\mathcal{I}_2}Y\to Bl_{\mathcal{I}_2}Bl_{\mathcal{I}_1}Y,
\end{align*}
from which we obtain the isomorphism $Bl_{\mathcal{I}_1\mathcal{I}_2}Y\cong Bl_{\mathcal{I}_2}Bl_{\mathcal{I}_1}Y$. The other isomorphism $Bl_{\mathcal{I}_1\mathcal{I}_2}Y\cong Bl_{\mathcal{I}_1}Bl_{\mathcal{I}_2}Y$ follows symmetrically.

For simplicity of notation, denote $Y_1=Bl_{\mathcal{I}_1}Y$,
$Y_2=Bl_{\mathcal{I}_2}Bl_{\mathcal{I}_1}Y$ and $Y_3=Bl_{\mathcal{I}_1\mathcal{I}_2}Y$.

\nid(i) We show the existence of $f$.
$$\xymatrix{&& Bl_{\mathcal{I}_1\mathcal{I}_2}Y \ar[d]^{\phi}\\
Bl_{\mathcal{I}_2}Bl_{\mathcal{I}_1}Y\ar@{.>}[urr]^f\ar[r]_-{\phi_2}&
Bl_{\mathcal{I}_1}Y\ar[r]_-{\phi_1} & Y}$$

By the universal property of blowing up, it suffices to show that  $$(\phi_1\phi_2)^{-1}(\I_1\I_2)\cdot
\O_{Y_2}$$ is an invertible sheaf. Indeed,
\begin{align*}
(\phi_1\phi_2)^{-1}(\I_1\I_2)\cdot
\O_{Y_2}&=\phi_2^{-1}\big{(}(\phi_1^{-1}\I_1\cdot\O_{Y_1})\cdot(\phi_1^{-1}\I_2\cdot\O_{Y_1})\big{)}\cdot
\O_{Y_2}\\
&=\big{(}\phi_2^{-1}(\phi_1^{-1}\I_1\cdot\O_{Y_1})\cdot \O_{Y_2}\big{)}\cdot\big{(}
\phi_2^{-1}(\phi_1^{-1}\I_2\cdot\O_{Y_1})\cdot \O_{Y_2}\big{)},
\end{align*}
in the last expression both factors are invertible sheaves, therefore the product is also invertible.

\nid(ii) We show the existence of $g$.
$$\xymatrix{Bl_{\mathcal{I}_1\mathcal{I}_2}Y \ar@{.>}[r]^-{g}\ar@{.>}[dr]^{h}\ar[ddr]_{\phi}&
Bl_{\mathcal{I}_2}Bl_{\mathcal{I}_1}Y\ar[d]^{\phi_2}\\
&Bl_{\mathcal{I}_1}Y\ar[d]^{\phi_1}\\ & Y}$$

Since $(\phi^{-1}\I_1\cdot\O_{Y_3})\cdot (\phi^{-1}\I_2\cdot\O_{Y_3})=
\phi^{-1}(\I_1\I_2)\cdot\O_{Y_3}$ is invertible, both $(\phi^{-1}\I_1\cdot\O_{Y_3})$ and
$(\phi^{-1}\I_2\cdot\O_{Y_3})$ are invertible. The invertibility of $(\phi^{-1}\I_1\cdot\O_{Y_3})$
implies the existence of $h$ by the universal property of blowing up. Then, since $\phi_2$ is the
blow-up of the ideal sheaf $(\phi_1^{-1}\I_2\cdot\O_{Y_1})$ and
$$h^{-1}(\phi_1^{-1}\I_2\cdot\O_{Y_1})\cdot \O_{Y_3}=\phi^{-1}\I_2\cdot\O_{Y_3}$$
is invertible, we can lift $h$ to $g$ by applying the universal property of blowing up again.
This completes the proof.
\end{proof}

Now we are ready to prove Theorem \ref{order} given in the Introduction.

\begin{proof}[Proof of Theorem \ref{order}]
  (i) We fix the indices of $\{G_i\}$ in an order that is compatible with inclusion relations, i.e., $i<j$ if
  $G_i\subset G_j$. Consider the blow-up $\phi: \widetilde{Y}:=Bl_{\I_1}Y\to Y$ where $\I_1$ is the ideal sheaf of $G_1$. Since $G_i$
  ($i>1$)
  either contains $G_1$ or is transversal to $G_1$ by Definition-Lemma \ref{fact 1}, the ideal sheaf $\phi^{-1}\I_{G_i}\cdot\mathcal{O}_{\widetilde{Y}}$ is either $\I_{\widetilde{G}_i}\cdot \I_E$ or
  $\I_{\widetilde{G}_i}$.
Since $\I_E$ is invertible, the blow-up of $\I_{\widetilde{G}_i}\cdot \I_E$ is isomorphic to the
blow-up of $\I_{\widetilde{G}_i}$, that is, the blow-up along the nonsingular subvariety
$\widetilde{G}_i$. By the same argument, each blow-up $Y_{k+1}\to Y_k$ is isomorphic to the blow-up
of the ideal sheaf $\psi^{-1}\I_{k+1}\cdot\O_{Y_k}$ where $\psi: Y_k\to Y$ is the natural morphism.
Therefore, by Lemma \ref{blowup ideal}, $$Y_\mathcal{G}\cong Bl_{\I_N}\cdots
Bl_{\I_2}Bl_{\I_1}Y\cong Bl_{\I_1\cdots \I_N}Y.$$

\nid(ii) Now assume the order of $\{G_i\}$ is not necessarily compatible with inclusion relations, but
satisfies $(*)$.

The proof is by induction with respect to $N$. The statement is
obviously true for $N=1$. Assume the statement (ii) is true for $N$. Consider
$\mathcal{G}_+=\mathcal{G}\cup \{G_{N+1}\}$. We need to show that
$Y_{\mathcal{G}_+}$ is isomorphic to the blow-up of $Y_\mathcal{G}$
along a nonsingular subvariety $\widetilde{G}_{N+1}$.

Suppose $F$ is minimal in $\mathcal{G}$, $\widetilde{Y}=Bl_FY$ and $\phi:
\widetilde{Y}\to Y$ is the natural morphism.
Proposition \ref{G+} implies that
$\widetilde{\mathcal{G}}_+:=\widetilde{\mathcal{G}}\cup\{\widetilde{G}_{N+1}\}$ is a
building set in $\widetilde{Y}$. There are two cases: if $F$ is not minimal in $\mathcal{G}_+$, then
$G_{N+1}$ must be minimal and $G_{N+1}\subsetneq F$, in this case
$\phi^{-1}\I_{G_{N+1}}\cdot\O_{\widetilde{Y}}=\I_{\widetilde{G}_{N+1}}$.
Now consider the case where $F$ is minimal in $\mathcal{G}_+$. In this case $G_{N+1}$
either contains $F$ or is transversal to $F$, thus
$\phi^{-1}\I_{G_{N+1}}\cdot\O_{\widetilde{Y}}$ is either
 $\I_{\widetilde{G}_{N+1}}\cdot \I_E$ or $\I_{\widetilde{G}_{N+1}}$.
In each situation, $\phi^{-1}\I_{G_{N+1}}\cdot\O_{\widetilde{Y}}$ is
isomorphic to $\I_{\widetilde{G}_{N+1}}$ up to an invertible sheaf.
Continue this procedure until all elements in $\mathcal{G}$ have
been blown up. Let $\psi: Y_\mathcal{G}\to Y$ be the natural
morphism. Then $\psi^{-1}\I_{G_{N+1}}\cdot\O_{Y_\mathcal{G}}$ is
isomorphic to the ideal sheaf of the nonsingular subvariety
$\widetilde{G}_{N+1}\subset Y_\mathcal{G}$ up to an invertible
sheaf, therefore the blow-up of $Y_\mathcal{G}$ along
$\psi^{-1}\I_{G_{N+1}}\cdot\O_{Y_\mathcal{G}}$ is isomorphic to the
blow-up of $Y_\mathcal{G}$ along $\widetilde{G}_{N+1}$. This
completes the proof.
\end{proof}

\ss {Examples of wonderful compactifications}

\sss{{Wonderful model of subspace arrangements}}\label{eg:DP} If we let
$Y=V$ be a finite dimensional vector space and let $\mathcal{S}$ be any finite collection of
subspaces of $V$ and construct the wonderful compactification of any building set of subspaces of $V$, one recovers the \emph{wonderful model of subspace
arrangements} by De Concini and Procesi.

It was discovered by De Concini and Procesi (\cite{DP})
that if a subset $\mathcal{G}\subseteq \mathcal{S}$ forms a so-called \emph{building set}, the closure of the natural locally
closed embedding
$$i: V\setminus \bigcup_{W\in \mathcal{G}} W\hookrightarrow V\tm \prod_{W\in \mathcal{G}}\P(V/W)$$
is a nonsingular variety birational to $V$. Moreover, the subspaces
in $\mathcal{S}$ are replaced by a normal crossing divisor.

\rmk This idea motivated a generalized definition of the so-called
wonderful conical compactifications for a complex manifold given by
MacPherson and Procesi (\cite{MP}). Our definition of wonderful
compactification is neither strictly general nor strictly less general than the wonderful compactification defined in \cite{MP}. On the one hand, our compactification does not include the conic case: all
the subvarieties involved in our paper are assumed to be
nonsingular. On the other hand, even over the complex field
$\mathbb{C}$, many arrangements of nonsingular varieties are not
conical. 

\sss{Fulton-MacPherson configuration spaces}\label{eg:FM} Let $X$ be a nonsingular variety, $Y=X^n$ and $\mathcal{G}$ be the set of diagonals of $X^n$. Our wonderful compactification gives the Fulton-MacPherson configuration space $X[n]$.

In \cite{FM}, Fulton
and MacPherson have constructed  a compactification
$X[n]$ of the configuration space $F(X, n)$ of $n$ distinct labeled
points in a nonsingular algebraic variety $X$. This compactification is related to
several areas of mathematics. In \cite{FM}, Fulton and
MacPherson use their compactification to construct a differential graded algebra which
is a model for the configuration space $F(X, n)$ in the sense of Sullivan.
Axelrod-Singer used an analogous construction in the setting of real
smooth manifolds in Chern-Simons perturbation theory (\cite{AS}).
Now we give a brief review of Fulton and
MacPherson's construction.

The configuration space $F(X, n)$ is an open subset of the Cartesian
product $X^n$ defined as the complement of all diagonals:
$$F(X,n):=X^n\setminus \bigcup_{|I|\ge 2}\Delta_I=\{(p_1,\dots,p_n)\in X^n| p_i\neq p_j, \forall\, i\neq j\}.$$

The construction of $X[n]$ by Fulton and MacPherson is inductive. They define $X[1]$ to be $X$ and $X[n+1]$ is the variety that results from a sequence of blow-ups of $X[n]\tm X$
along nonsingular subvarieties corresponding to all diagonals $\Delta_I$ where $I\subseteq[n+1]$, $|I|\ge 2$ and $I$ contains the number $n+1$.

For example, $X[2]$ is the blow-up of $X^2$ along the diagonal
$\Delta_{12}$. The variety $X[3]$ is obtained from a sequence of blow-ups of $X[2]\tm X$ along non-singular subvarieties
corresponding to $\{\Delta_{123}; \Delta_{13}, \Delta_{23}\}$. More specifically, denote by $\pi$
the blow-up $X[2]\tm X\to X^3$, we blow up first along $\pi^{-1}(\Delta_{123})$, then along the
strict transforms of $\Delta_{13}$ and $\Delta_{23}$ (the two strict transforms are disjoint, so
they can be blown up in any order). In general, the order of blow-ups in the construction of $X[n]$ can be expressed as  $$\Delta_{12}, \Delta_{123}, \Delta_{13}, \Delta_{23}, \Delta_{1234}, \Delta_{124}, \Delta_{134}, \Delta_{234}, \Delta_{14},\Delta_{24},\Delta_{34},\Delta_{12345},\Delta_{1235},\dots $$
It is easy to verify that the above sequence satisfies $(*)$ in Theorem \ref{order}, therefore the resulting variety $X[n]$ is indeed the wonderful compactification $Y_\mathcal{G}$.
Theorem \ref{order} also implies that $X[n]$ can be obtained from a more symmetric sequence
of blow-ups in the order of ascending dimension:
$$\Delta_{12\cdots n}, \Delta_{12\cdots (n-1)},\dots, \Delta_{23\cdots n},\dots,\Delta_{12},\dots, \Delta_{(n-1),n}.$$
This more symmetric order of blow-ups is given by De Concini and
Procesi \cite{DP}, MacPherson and Procesi \cite{MP}, and Thurston
\cite{Thurston}.

In fact, graphs can be used to clarify the condition $(*)$ by using
the Kuperberg-Thurston's compactification (cf. the discussion after
Proposition \ref{graph} below).

\sss{Kuperberg-Thurston's compactifications}\label{eg:KT} In their
paper \cite{KT}, Kuperberg and Thurston have constructed an
interesting compactification of the configuration space $F(X,n)$. Their
construction is for real smooth manifolds and we adapt their
compactification in this section to a nonsingular algebraic variety.

Let $\Gamma$ be a (not necessarily connected) graph with $n$ labeled
vertices such that  $\Gamma$ has no self-loops and multiple edges.
Denote by $\Delta_{\Gamma}$
 the polydiagonal in $X^n$ where $x_i=x_j$ if $i, j$ are connected in $\Gamma$. We call a graph $\Gamma$
\emph{vertex-2-connected} if the graph is connected and will still
be connected if we remove any vertex. In particular, a single edge
is vertex-2-connected.

In \cite{KT} the authors state and sketch a proof that blowing up
along $\Delta_{\Gamma'}$ for all vertex-2-connected subgraphs
$\Gamma'\subseteq\Gamma$ gives a compactification $X^\Gamma$. When
$\Gamma$ is the complete graph with $n$ vertices (i.e. any two vertices
are joined with an edge), the compactification $X^\Gamma$ is exactly
the Fulton-MacPherson compactication $X[n]$.

Kuperburg-Thurston's compactification $X^\Gamma$ is a
special case of the wonderful compactification of an arrangement of
subvarieties given in this paper. Indeed, let $Y=X^n$, let
$$\mathcal{G}:=\{\Delta_{\Gamma'}: \Gamma'\subseteq\Gamma \hbox{ is vertex-2-connected}\},$$ and let
$\mathcal{S}$ be the set of polydiagonals of $X^n$ obtained by intersecting only the diagonals in
$\mathcal{G}$.

\begin{prop}
  In the above notation, $\mathcal{G}$ is a building set with respect to $\mathcal{S}$. Therefore Kuperburg-Thurston's compactification is the wonderful compactification $Y_{\mathcal{G}}$.
\end{prop}

\begin{proof}
The proof is in two steps.

\nid(i) We call $\Gamma'\subseteq\Gamma$ a \emph{full subgraph} if the following are satisfied.
\begin{itemize}
\item $\Gamma'$ contains all vertices in $\Gamma$.
\item For any edge $e\in \Gamma$, if its endpoints $p$ and $q$ are in the same connected component of $\Gamma'$, then $e\in\Gamma'$.
\end{itemize}
Then there is a one-one correspondence between the set of all full
subgraphs of $\Gamma$ and the set $\mathcal{S}$. The correspondence
is given by mapping a full subgraph $\Gamma'$ to $\Delta_{\Gamma'}$.

\nid(ii) Any full subgraph $\Gamma'$ has a unique decomposition into vertex-2-connected subgraphs
$\Gamma_1,\dots,\Gamma_k$. Notice that
$\Delta_{\Gamma_1},\dots,\Delta_{\Gamma_k}$ are the minimal elements
in $\mathcal{G}$ which contain $\Delta_{\Gamma'}$, and they
intersect transversally with the intersection $\Delta_{\Gamma'}$.
Therefore $\mathcal{G}$ is a building set by Definition \ref{def
building set subvariety}.
\end{proof}

\rmk It is also easy to describe a $\mathcal{G}$-nest: it
corresponds to a set of vertex-2-connected subgraphs of $\Gamma$,
where for any two subgraphs $\Gamma_1$ and $\Gamma_2$ one of the following holds.

\begin{itemize}
 \item[(i)] $\Gamma_1$ and $\Gamma_2$ are disjoint, or
 \item[(ii)] $\Gamma_1$ and $\Gamma_2$ intersect at one vertex, or
 \item[(iii)] $\Gamma_1\subseteq\Gamma_2$ or $\Gamma_2\subseteq \Gamma_1$.
\end{itemize}

\medskip

The following proposition describes the relation between $X^{\Gamma_1}$ and $X^{\Gamma_2}$ for  $\Gamma_1\subsetneq\Gamma_2$, which will help us to understand the construction of Fulton and MacPherson configuration spaces.
\begin{prop}\label{graph}
  Let $\Gamma_1\subsetneq\Gamma_2$ be two (not necessarily connected) graphs with $n$ labeled vertices without
  self-loops and multiple edges. Then $X^{\Gamma_2}$ can be obtained by a sequence of blow-ups of $X^{\Gamma_1}$
  along nonsingular centers.

  One such order is given as follows. Let $\{\Gamma''_{j}\}_{j=1}^t$ be the set of all vertex-2-connected
  subgraphs of $\Gamma_2$ that are not contained in $\Gamma_1$. Arrange the index such that $i<j$ if the
  number of vertices of $\Gamma''_{i}$ is greater than the number of vertices of $\Gamma''_{j}$. Then $X^{\Gamma_2}$ can be obtained by blowing up along the following nonsingular centers   $\widetilde{\Delta}_{\Gamma''_{1}},\widetilde{\Delta}_{\Gamma''_{2}},\dots,\widetilde{\Delta}_{\Gamma''_{t}}$.

\end{prop}

\begin{proof}
  Let $\{\Gamma'_i\}_{i=1}^s$ be the set of all vertex-2-connected
  subgraphs of $\Gamma_1$, arrange the indices such that $i<j$ if the
  number of vertices of $\Gamma'_i$ is greater than the number of vertices of $\Gamma'_j$.

  It is easy to verify that $\{\Delta_{\Gamma'_1},\dots, \Delta_{\Gamma'_s}, \Delta_{\Gamma''_1},\dots,\Delta_{\Gamma''_t} \}$
  satisfies $(*)$ in Theorem \ref{order}. Apply Theorem \ref{order}, we know that $X^{\Gamma_2}$ is
  the blowup of $X^n$ along nonsingular centers
  $$\Delta_{\Gamma'_1}, \widetilde{\Delta}_{\Gamma'_2}\dots,  \widetilde{\Delta}_{\Gamma'_s}, \widetilde{\Delta}_{\Gamma''_1},\dots,\widetilde{\Delta}_{\Gamma''_t}.$$
  On the other hand, after the first $s$ blow-ups, we get $X^{\Gamma_1}$. Therefore
  $X^{\Gamma_2}$ can be obtained by a sequence of blow-ups along nonsingular centers $\widetilde{\Delta}_{\Gamma''_1},\dots,\widetilde{\Delta}_{\Gamma''_t}$.
  This completes the
  proof.
\end{proof}

With the above proposition, Fulton and MacPherson's original construction of $X[n]$ can be
understood as specifying a chain of graphs:
\begin{figure}[h]
\hskip10mm\xy
    (88,-8)*{\mycaption{7in}{
            Fulton and MacPherson's construction of $X[4]$.}};
    (5,0)="m";
    "m"+(5,0)  *{X^4};
    "m"+(0,10)  *{\bullet};
    "m"+(0,20)  *{\bullet};
    "m"+(10,10)  *{\bullet};
    "m"+(10,20)  *{\bullet};
    "m"+(0,23)  *{\scriptstyle 1};
    "m"+(0,7)  *{\scriptstyle 2};
    "m"+(10,23)  *{\scriptstyle 4};
    "m"+(10,7)  *{\scriptstyle 3};
    "m"+(20,15)  *{\longrightarrow};
(35,0)="m";
    "m"+(5,0)  *{X[2]\times X^2};
    "m"+(0,10)="1"  *{\bullet};
    "m"+(0,20)="2"  *{\bullet};
    "m"+(10,10)="3"  *{\bullet};
    "m"+(10,20)="4"  *{\bullet};
    "m"+(0,23)  *{\scriptstyle 1};
    "m"+(0,7)  *{\scriptstyle 2};
    "m"+(10,23)  *{\scriptstyle 4};
    "m"+(10,7)  *{\scriptstyle 3};
    "2";         "1"**\dir{-};
    "m"+(20,15)  *{\longrightarrow};
(65,0)="m";
    "m"+(5,0)  *{X[3]\times X};
    "m"+(0,10)="1"  *{\bullet};
    "m"+(0,20)="2"  *{\bullet};
    "m"+(10,10)="3"  *{\bullet};
    "m"+(10,20)="4"  *{\bullet};
    "m"+(0,23)  *{\scriptstyle 1};
    "m"+(0,7)  *{\scriptstyle 2};
    "m"+(10,23)  *{\scriptstyle 4};
    "m"+(10,7)  *{\scriptstyle 3};
    "2";         "1"**\dir{-};
    "3";         "1"**\dir{-};
    "3";         "2"**\dir{-};
    "m"+(20,15)  *{\longrightarrow};
(95,0)="m";
    "m"+(5,0)  *{X[4]};
    "m"+(0,10)="1"  *{\bullet};
    "m"+(0,20)="2"  *{\bullet};
    "m"+(10,10)="3"  *{\bullet};
    "m"+(10,20)="4"  *{\bullet};
    "m"+(0,23)  *{\scriptstyle 1};
    "m"+(0,7)  *{\scriptstyle 2};
    "m"+(10,23)  *{\scriptstyle 4};
    "m"+(10,7)  *{\scriptstyle 3};
    "2";         "1"**\dir{-};
    "3";         "1"**\dir{-};
    "3";         "2"**\dir{-};
    "4";         "1"**\dir{-};
    "4";         "2"**\dir{-};
    "4";         "3"**\dir{-};
\endxy
\end{figure}
indeed, by Proposition \ref{graph}, the first arrow corresponds to blowing up $X[4]$ along $\Delta_{12}$, the second
corresponds to blowing up along $\widetilde{\Delta}_{123}$,
$\widetilde{\Delta}_{13}$ and $\widetilde{\Delta}_{23}$ (which correspond to all the vertex-2-connected subgraphs which are not in the previous graph), the last arrow corresponds to
blowing up along  $\widetilde{\Delta}_{1234}$, $\widetilde{\Delta}_{124}$,
$\widetilde{\Delta}_{134}$, $\widetilde{\Delta}_{234}$,
$\widetilde{\Delta}_{14}$, $\widetilde{\Delta}_{24}$,
$\widetilde{\Delta}_{34}$.

On the other hand, the symmetric construction of $X[4]$ corresponds
to the chain containing only two graphs: the first graph and last
graph in Figure 1.
\medskip

To illustrate the idea a little more, we construct $X[3]$
corresponding to the chain of graphs in the figure below. The first step is to blow up along $\Delta_{12}$, the second step is to blow up along 
$\widetilde{\Delta}_{23}$, the final step is to blow up along $\widetilde{\Delta}_{123}$ and $\widetilde{\Delta}_{13}$. Each
blow-up is along a nonsingular subvariety.

\begin{figure}[h]
\hskip10mm\xy
    (88,-8)*{\mycaption{7in}{
            A new construction of $X[3]$.}};
    (5,0)="m";
    "m"+(5,0)  *{X^3};
    "m"+(0,10)  *{\bullet};
    "m"+(0,20)  *{\bullet};
    "m"+(10,10)  *{\bullet};
    "m"+(0,23)  *{\scriptstyle 1};
    "m"+(0,7)  *{\scriptstyle 2};
    "m"+(10,7)  *{\scriptstyle 3};
    "m"+(20,15)  *{\longrightarrow};
(35,0)="m";
    "m"+(5,0)  *{X[2]\times X};
    "m"+(0,10)="1"  *{\bullet};
    "m"+(0,20)="2"  *{\bullet};
    "m"+(10,10)="3"  *{\bullet};
    "m"+(0,23)  *{\scriptstyle 1};
    "m"+(0,7)  *{\scriptstyle 2};
    "m"+(10,7)  *{\scriptstyle 3};
    "2";         "1"**\dir{-};
    "m"+(20,15)  *{\longrightarrow};
(65,0)="m";
    "m"+(5,0)  *{Bl_{\widetilde{\Delta}_{23}}(X[2]\times X)};
    "m"+(0,10)="1"  *{\bullet};
    "m"+(0,20)="2"  *{\bullet};
    "m"+(10,10)="3"  *{\bullet};
    "m"+(0,23)  *{\scriptstyle 1};
    "m"+(0,7)  *{\scriptstyle 2};
    "m"+(10,7)  *{\scriptstyle 3};
    "2";         "1"**\dir{-};
    "3";         "1"**\dir{-};
    "m"+(20,15)  *{\longrightarrow};
(95,0)="m";
    "m"+(5,0)  *{X[3]};
    "m"+(0,10)="1"  *{\bullet};
    "m"+(0,20)="2"  *{\bullet};
    "m"+(10,10)="3"  *{\bullet};
    "m"+(0,23)  *{\scriptstyle 1};
    "m"+(0,7)  *{\scriptstyle 2};
    "m"+(10,7)  *{\scriptstyle 3};
    "2";         "1"**\dir{-};
    "3";         "1"**\dir{-};
    "3";         "2"**\dir{-};
\endxy
\end{figure}

\sss{Moduli space $\overline{M}_{0,n}$ of rational curves with $n$ marked points
}\label{eg:Moduli}

The moduli space $\overline{M}_{0,n}$ is the wonderful compactification of
$((\mathbb{P}^1)^{n-3},\mathcal{S},\mathcal{G})$ where $\mathcal{G}$
is set of all diagonals and augmented diagonals
$$\Delta_{I,a}:=\{(p_4,\cdots,p_n)\in (\mathbb{P}^1)^{n-3}: p_i=a, \forall i\in I\}$$ for
$I\subseteq \{4,\dots,n\}$, $|I|\ge 2$ and $a\in\{0, 1,\infty\}$.
$\mathcal{S}$ is the set of all intersections of elements in
$\mathcal{G}$.

This is a immediate consequence of Theorem \ref{order} applied to
Keel's construction \cite{Keel1}. Indeed, Keel gives the
construction of $\overline{M}_{0,n}$ by a sequence of blow-ups in
the following order:

$\Delta_{45,0},\Delta_{45,1},\Delta_{45,\infty},\Delta_{456,0},\Delta_{456,1},\Delta_{456,\infty},\dots,\Delta_{46,0},\dots,\Delta_{456}, \dots$

\nid To be more precise: for $I$ such that $\max I=5$, blow up along $\Delta_{I,a}$ for those $I$ such that $|I|=2$; for $\max I=6$, blow up $\Delta_{I,a}$ for $|I|=3$, then $\Delta_I$ for $|I|=3$; in general, for $\max I=k$, blow up $\Delta_{I,a}$ for $|I|=k-3$, then  $\Delta_{I,a}$ for $|I|=k-4$ and $\Delta_I$ for $|I|=n-3$, then  $\Delta_{I,a}$ for $|I|=k-5$ and $\Delta_I$ for $|I|=n-4$, ...

It is easy to check the order satisfies $(*)$
in Theorem \ref{order}. So $\overline{M}_{0,n}$ is a wonderful compactification.

Notice that the above diagonals and augmented diagonals in
$\P^{n-3}$ are just the restrictions of diagonals in $(\P^1)^n$ to
the codimension three subvariety $$Y=\{(p_1, p_2,\cdots, p_n)\in
(\P^1)^n: p_1=0, p_2=1, p_3=\infty\}.$$

Now blow up all the diagonals of $(\P^1)^n$ in order of increasing dimension and compare with the construction of Fulton-MacPherson configuration space, we get a
relation between $\overline{M}_{0,n}$ and the Fulton-MacPherson
space $\P^1[n]$ in Corollary \ref{corM0n}.

\sss{Ulyanov's compactifications}\label{eg:U} Closely related to
Fulton and MacPherson's compactification, Ulyanov has discovered
another compactification of the configuration space $F(X,n)$, which
he denoted by $X\pd<n>$ (\cite{Ulyanov}). The construction consists
of blowing up more subvarieties in $X^n$ than Fulton-MacPherson's
construction does, i.e., it blows up not only diagonals but also polydiagonals. 
The order of the blow-ups in \cite{Ulyanov} is the ascending order
of the dimension. For example, $X\pd<4>$ is the
blow-up of $X^4$ along polydiagonals in the following order:
$$(1234),(123),(124),(134),(234),(12,34),(13,24),(14,23),(12),\dots,(34).$$

The polydiagonal compactification  $X\pd<n>$ shares many similar properties
with Fulton-MacPherson's compactification. However, one difference is that in the case
of characteristic 0, the
isotropy  group of any point in $X\pd<n>$ is abelian under the symmetric group action, while the
isotropy group of a point in $X[n]$ is not necessarily abelian (but always solvable) under the symmetric group action.

\medskip
\sss{Hu's compactifications v.s. minimal
compactifications}\label{eg:H} We now consider the general situation
where $Y$ is nonsingular with an arrangement of subvarieties
$\mathcal{S}$. By blowing up all $S\in \mathcal{S}$ in order of
ascending dimesion, we get a nonsingular variety $Bl_\mathcal{S}Y$
(\cite{Hu}). Define $Y^\circ:=Y\setminus\cup_{S\in\mathcal{S}} S$,
the open stratum of $Y$. It is isomorphic to an open subset of
$Bl_\mathcal{S}Y$. Hu showed \emph{\begin{enumerate}
\item[(i)] The boundary $Bl_\mathcal{S}Y\setminus Y^\circ=\cup_{S\in\mathcal{S}} D_S$ is a simple normal crossing divisor.
\item[(ii)] For any $S_1,..., S_n\in\mathcal{S}$, the intersection of $D_{S_1}\dots D_{S_k}$ is nonempty if and only if $\{S_i\}$ forms a chain,
i.e., $S_1\subseteq\dots\subseteq S_k$ with a rearrangement of
indices if necessary.
\end{enumerate}}
\medskip

Hu's compactification generalized Ulyanov's polydiagonal
compactification and is a special case of the wonderful
compactification of arrangement of subvarieties given in this paper
where the building set $\mathcal{G}=\mathcal{S}$. (In this special
case, a $\mathcal{G}$-nest is simply a chain of subvarieties.)

Fixing an arrangement $\mathcal{S}$, Hu's compactification
$Y_\mathcal{S}$ is the maximal wonderful compactification. Indeed, it is not hard to show that for any building set $\mathcal{G}$ of $\mathcal{S}$, the natural birational map $Y_\mathcal{S}\to Y_\mathcal{G}$ is a morphism. On the other extreme, there exists a
minimal wonderful compactification for $\mathcal{S}$, which can be defined by the set of so-called irreducible elements in $\mathcal{S}$.

\begin{df}
  An element $G$ in $\mathcal{S}$ is called {\em reducible} if
  there are $G_1$, $\dots$, $G_k \in \mathcal{S}$ ($k\ge 2$) with $G=G_1\pitchfork\dots\pitchfork
  G_k$ and for every $G'\supseteq G$ in $\mathcal{S}$, there exist
  $G'_i\in\mathcal{S}$, $G'_i\supseteq G_i$ for $1\le i\le k$, such that $G'=G'_1\pitchfork\dots\pitchfork
  G'_k$.

  $G\in\mathcal{S}$ is called {\em irreducible} if it is not reducible.
\end{df}

By the same method as in \cite{DP}, we can show that the irreducible
elements in $\mathcal{S}$ form a building set, denoted by $\mathcal{G}_{min}$,
and that every building set $\mathcal{G}$ of the arrangement
$\mathcal{S}$ contains $\mathcal{G}_{min}$. It is not hard to show that the natural birational map $Y_\mathcal{G}\to
Y_{\mathcal{G}_{min}}$ is a morphism. 

Out of the previous examples, Fulton-MacPherson
configuration spaces, Kuperberg-Thurston's compactifications and the moduli
space $\overline{M}_{0,n}$ are minimal wonderful compactifications. Ulyanov's polydiagonal compactifications are maximal.

 \ss{Appendix}

\sss{Clean intersection v.s. transversal intersection}\label{intersections}

Let $Y$ be a nonsingular variety. For a nonsingular subvariety $A$ (more generally, a subscheme whose connected components are nonsingular
subvarieties) of $Y$, denote by $T_A$ the total space of the tangent bundle of $A$, denote by $T_{A,y}$ the tangent space of $A$ at the point $y\in A$. For a point $y\notin A$, define $T_{A,y}$ to be $T_y$, the tangent space of $Y$ at $y$. (This artificial definition will simplify the definition of transversal intersection.) In this paper, $T_{A,y}$ will be seen as a subspace of $T_y$ and $T_A$ will be seen as a subvariety of $T_Y$.

\ssss{Clean intersection.} The notion of cleanness can be traced back to Bott \cite{Bott} in the setting of differential geometry.

We say that the intersection of two nonsingular subvarieties $A$ and $B$ is \emph{clean} if the set-theoretic intersection $A\cap B$ is a nonsingular
subvariety (or, more generally, a scheme whose connected components are nonsingular subvarieties)  and satisfies the condition
$$T_{A\cap B,y}=T_{A,y}\cap T_{B,y}, \quad \forall y\in A\cap B.$$

The following lemma gives a useful criterion for the cleanness of intersections.

\begin{lem}\label{ideal_clean}
  Suppose $A$ and $B$ are nonsingular closed subvarieties of $Y$ and the intersection $C=A\cap B$ is a disjoint union of nonsingular subvarieties.
  Let $\I_A$ (resp. $\I_B$, $\I_C$) denote the ideal sheaf of $A$ (resp. $B$, $C$). Then the following are equivalent:
\begin{itemize}
\item[(i)] The subvarieties $A$ and $B$ intersect cleanly.

\item[(ii)]  $\I_A+\I_B=\I_C$.
\end{itemize}
In other words, two subvarieties intersect cleanly if and only if their scheme-theoretic intersection is nonsingular.
\end{lem}

\begin{proof}
Condition (i) is equivalent to
\begin{equation}\label{e2}
T_{A,y}\cap T_{B,y}=T_{C,y}, \quad \forall y\in A\cap B.
\end{equation}
By definition of tangent space,
$$T_{A,y}=\{v\in T_y|\, df(v)=0, \forall f\in (\I_A)_y\}.$$
Define $\phi: \mathfrak{m}_y\to \mathfrak{m}_y/\mathfrak{m}_y^2$ to be the natural quotient. Then $T_{A,y}=\phi\big{(}(\I_A)_y\big{)}^\perp$, the
annihilator of $\phi\big{(}(\I_A)_y\big{)}$ in the dual space $(\mathfrak{m}_y/\mathfrak{m}_y^2)^*\cong T_y$. Therefore condition (\ref{e2}) is
equivalent to
$$\phi\big{(}(\I_A)_y\big{)}^\perp\cap\phi\big{(}(\I_B)_y\big{)}^\perp=\phi\big{(}(\I_C)_y\big{)}^\perp, \quad \forall y\in A\cap B,$$
which is equivalent to
$$\phi\big{(}(\I_A)_y\big{)}+\phi\big{(}(\I_B)_y\big{)}=\phi\big{(}(\I_C)_y\big{)}, \quad \forall y\in A\cap B. $$
Since $\phi\big{(}(\I_A)_y\big{)}=((\I_A)_y+\mathfrak{m}_y^2)/\mathfrak{m}_y^2$ (similarly for $B$ and $C$),  the above condition is equivalent to
\begin{equation}\label{e3}
(\I_A)_y+(\I_B)_y+\mathfrak{m}_y^2=(\I_C)_y+\mathfrak{m}_y^2, \hbox{\quad $\forall y\in A\cap B$.}
\end{equation}
On the other hand, two ideal sheaves on $Y$ are the same if and only if their germs coincide at every closed point $y\in Y$. So condition (ii) is
equivalent to
\begin{equation}\label{e1} (\I_A)_y+(\I_B)_y=(\I_C)_y, \hbox{\quad $\forall y\in Y$}
\end{equation}
where $\I_y$ denote the germ of a sheaf $\I$ at point $y$. Therefore it suffices to show that (\ref{e3})$\Leftrightarrow$(\ref{e1}).

Obviously (\ref{e1}) $\Rightarrow$ (\ref{e3}). To see the implication (\ref{e3}) $\Rightarrow$ (\ref{e1}), observe first that condition (\ref{e1})
holds for $y\notin A\cap B$ and the inclusion $(\I_A)_y+(\I_B)_y\subseteq(\I_C)_y$ holds for $y\in A\cap B$. Thus it remains to show
$(\I_A)_y+(\I_B)_y\supseteq(\I_C)_y$ holds for $y\in A\cap B$. Using local parameters it can be checked that $(\I_C)_y\cap
\mathfrak{m}_y^2=(\I_C)_y\mathfrak{m}_y$. Condition (\ref{e3}) then implies the surjection
$$(\I_A)_y+(\I_B)_y\twoheadrightarrow
\big{(}(\I_C)_y+\mathfrak{m}_y^2\big{)}/\mathfrak{m}_y^2\stackrel{\cong}{\to}(\I_C)_y/\big{(}(\I_C)_y\cap\mathfrak{m}_y^2\big{)}\stackrel{\cong}{\to}
(\I_C)_y/(\I_C)_y\mathfrak{m}_y.
$$
Hence $(\I_A)_y+(\I_B)_y+(\I_C)_y\mathfrak{m}_y=(\I_C)_y$. Apply Nakayama's lemma, we have $(\I_A)_y+(\I_B)_y=(\I_C)_y$. This completes the proof.
\end{proof}

\ssss{Transversal intersection.} By definition, $A$
and $B$ intersect transversally (denoted by $A\pitchfork B$) if
$T_{A,y}^\perp+ T_{B,y}^\perp$ form a direct sum in the dual
space $T_y^*\cong\mathfrak{m}_y/\mathfrak{m}_y^2$ of $T_y$ for any point $y\in Y$,  or equivalently, if
$$T_y=T_{A,y}+T_{B,y}, \quad \forall y\in Y.$$

More generally, we define that a finite collection of $k$ nonsingular
subvarieties $A_1$, $\dots$, $A_k$  intersect
transversally (denoted by $A_1\pitchfork
A_2\pitchfork\cdots\pitchfork A_k$) if either $k=1$ or for any $y\in
Y$,
$$T_{A_1,y}^\perp+ T_{A_2,y}^\perp+\cdots+T_{A_k,y}^\perp$$  form a direct sum in
$T_y^*$; or equivalently, if
$${\rm codim}(\bigcap_{i=1}^k T_{A_i,y}, T_y)=\sum_{i=1}^k{\rm codim}(A_i, Y);$$
or equivalently, if for any $y\in Y$, there exists a system of local
parameters $x_1,\dots,x_n$ on $Y$ at $y$ which are regular on an
affine neighborhood $U$ of $y$ such that $y$ is defined by the maximal
ideal $(x_1,\dots, x_n)$, and there exist integers $0=r_0\le
r_1\le\cdots\le r_k\le n$ such that the subvariety $A_i$ is defined by the ideal
$$(x_{r_{i-1}+1}, x_{r_{i-1}+2}, \cdots, x_{r_{i}}), \quad \forall 1\le i\le k.$$
(In case that $r_{i-1}=r_i$, the ideal is assumed to be the ideal
containing units, which means geometrically that the restriction of
$A_i$ to $U$ is empty.)

\ssss{Transversal intersection $\Rightarrow$ clean intersection.} If $A$
and $B$ intersect transversally, we can choose local parameters
around any point $y\in A\cap B$ such that $y$ is the origin and the restriction of $A$
and $B$ are defined by local parameters. Then it is obvious that
$T_{A\cap B,y}=T_{A,y}\cap T_{B,y}$, $\forall y\in A\cap B$.

\ssss{Transversal intersection at one point + clean intersection $\Rightarrow$ transversal intersection.}

\begin{lem}\label{appendix lemma}
  Let $A_1$ and $A_2$ be two nonsingular closed subvarieties of $Y$ that intersect cleanly along a closed nonsingular subvariety $A$. If $A_1$
  and $A_2$ intersect transversally at a point $y_0\in A$, then they intersect transversally (at every point $y\in A$).

In general, let $A_1,\dots, A_k$ be subvarieties in a simple
arrangement $\mathcal{S}$ (cf. Definition \ref{def arrangement}), let $A=\cap_{i=1}^k{A_i}$. If $A_1,\dots,A_k$ intersect transversally at a point
$y_0\in A$, then they intersect transversally (at every point).
\end{lem}

\begin{proof}
 We prove the general case. Without loss of
generality, we need only to prove the transversality for points in $A$. The
irreducibility of $A_i$ and $A$ implies $\dim T^\perp_{A_i,y}=\dim
T^\perp_{A_i,y_0}$ and $\dim T^\perp_{A,y}=\dim T^\perp_{A,y_0}$. By
the definition of clean intersection, we have
$$T^\perp_{A_1,y}+\cdots+ T^\perp_{A_k,y}=T^\perp_{A,y}.$$
On the other hand, by the transversality condition at point $y_0$,
$$T^\perp_{A_1,y_0}\oplus\cdots\oplus
T^\perp_{A_k,y_0}=T^\perp_{A,y_0}.$$ Comparing the dimensions of the
above two equalities,  the left hand side of the first equality must
form a direct sum, therefore $A_1,\dots,A_k$ intersect transversally
at $y$.
\end{proof}

\ssss{Examples/nonexamples of clean and transversal intersections.}
\begin{itemize}
  \item $k (\le n)$ hyperplanes $H_i$ in $\mathbb{A}^n$ defined by $x_i=0$ intersect transversally,
  therefore any two of them intersect cleanly.
  \item Two (not necessarily distinct) lines in $\mathbb{A}^3$ passing through the origin intersect
  cleanly but not transversally.
  \item In $\mathbb{A}^2$, the intersection of the parabola $y=x^2$ and the line $y=0$ is not
  clean, therefore not transversal.
\end{itemize}

\sss{Proofs of statements in previous sections}\label{appendix2}

\begin{proof}[Proof of Lemma \ref{fact 2}]

  It is convenient to carry out the proof using the cotangent space $T_y^*$. We
  use the same notation $\mathcal{G}^*_y$, $S$, $T_{S,y}^\perp$, $T_i^\perp$ as in the remark after Definition \ref{def building set
  subvariety}. By \cite{DP} \S2.3 Theorem (2), the definition of building set
implies the following:

If $S'\in \mathcal{S}$ is such that $S'\supseteq S$, then
$$T_{S',y}^\perp=\bigoplus_{i=1}^k\big{(}T_{S',y}^\perp\cap T_i^\perp\big{)},$$ moreover, if
$T_{S',y}^\perp=T_1^{'\perp}\oplus\cdots \oplus T_s^{'\perp}$ where
$T_1^{'\perp},\dots, T_s^{'\perp}$ are the maximal elements in
$\mathcal{G}^*_y$ contained in $T_{S',y}^\perp$, then each term
$(T_{S',y}^\perp\cap T_i^\perp)$ is a direct sum of some
$T_j^{'\perp}$'s.

 Fix a point $y\in S$. To show (i) it is enough to show the following:

  Suppose that $T_1^\perp,\dots, T_k^\perp$ are all the maximal elements in $\mathcal{G}^*_y$ which are contained in
  $T_{S,y}^\perp$. Suppose  $m$ is an integer such that $1\le m\le k$,  and define  $T^\perp:=T_1^\perp\oplus\cdots\oplus T_m^\perp$. Then $T_1^\perp,\dots, T_m^\perp$
  are all the maximal elements in $\mathcal{G}^*_y$ which are contained in $T^\perp$.

To show (ii) it is equivalent to show:

  Suppose $T_1^\perp,\dots, T_k^\perp$ are all the maximal elements in $\mathcal{G}^*_y$ that are
 contained in $T_{S,y}^\perp$. Suppose
  $T^\perp\in\mathcal{G}^*_y$ is maximal, $T^\perp\supseteq T_1^\perp,\dots,T_m^\perp$ and
  $T^\perp\nsupseteq T_{m+1}^\perp, \dots, T_k^\perp$. Then $T^\perp, T_{m+1}^\perp,\dots, T_k^\perp$ are all the
  maximal elements in $\mathcal{G}^*_y$ which are contained in $T^\perp+ T_{S,y}^\perp$.


Both statements can be shown by routine linear algebra.
\end{proof}

\begin{proof}[Proof of Definition-Lemma \ref{fact 1}]
\nid(i) This part follows directly from the definition of building set: if $F$ is disjoint from $G$ then
of course $F\pitchfork G$; otherwise $G$ contains some $\mathcal{G}$-factor of $F\cap G$. But a
$\mathcal{G}$-factor of $F\cap G$ is either $F$ or is transversal to $F$ (which implies
$G\pitchfork F$).

\nid(iii) Define $A=\cap_{i=1}^mG_i$ and $B=\cap_{i=m+1}^kG_i$. We claim that
$A\supseteq F$ and $B\pitchfork F$. $A=\cap_{i=1}^mG_i\supseteq F$
is because of the definition of $m$. Lemma \ref{fact 2} (ii) asserts
$F\pitchfork G_{m+1}\pitchfork\cdots\pitchfork G_k$, so $F$ is
transversal to $B$. Then (iii) follows from Lemma \ref{fact 2} (i).

\nid(ii) The above proof of (iii) shows the existence of an $F$-factorization. Now we show that such an factorization is unique. Given another factorization $S=A'\cap B'$ such that $A'\supseteq
F$ and $B'\pitchfork F$. Since $B'\supseteq F\cap B'=F\cap S$ and the $\mathcal{G}$-factors of
$F\cap S$ are $F, G_{m+1},\dots, G_k$ by Lemma \ref{fact 2} (ii), each $\mathcal{G}$-factor $G'$
of $B'$ contains $F$ or $G_i$ for some $m+1\le i\le k$. But $B'\pitchfork F$ implies $G'\pitchfork
F$, hence $G'\nsupseteq F$. So $G'\supseteq G_i$ for some $m+1\le i\le k$. Intersecting
all the $\mathcal{G}$-factors $G'$
of $B'$, we have $B'=\cap G'\supseteq \cap_{i=m+1}^kG_i=B.$ Fixing a point $y\in F\cap S$, we have
$$T_{F,y}^\perp\oplus T_{B,y}^\perp=T_{F,y}^\perp\oplus T_{B',y}^\perp$$
and  $T_{B,y}^\perp\supseteq T_{B',y}^\perp$, therefore
$T_{B,y}^\perp=T_{B',y}^\perp$ hence $B=B'$. Similarly $A=A'$.

\nid(iv) Suppose the $F$-factorization of $S\cap S'$ is $A''\cap B''$. Then $F\cap B''$ is the
$F$-factorization of the intersection. Since $B\supseteq (F\cap S)=(F\cap B'')$ but $B\pitchfork
F$, $B\supseteq B''$. Similarly $B'\supseteq B''$. So $B\cap B'\supseteq B''$. By an analogous
argument using the dual of the tangent space as in the proof of (ii), we can show that $B\cap B'=B''$. So $F\pitchfork (B\cap B')$. Then it is easy to see that $A''=A\cap A'$ and the $F$-factorization of $S\cap S'$ is indeed $(A\cap A')\cap (B\cap B')$.
\end{proof}

\begin{proof}[Proof of Lemma \ref{basic}]
We give only the proof of (iii) since (ii) and (iv) can be proved
similarly, while
 (i), (v) and (vi) can be easily checked using a system of local
 parameters.

In the complement of the exceptional divisor $E$, we have
$$(\widetilde{A}_1\cap\widetilde{A}_2)\setminus E\cong(A_1\setminus F)\cap (A_2\setminus F)=(A_1\cap A_2)\setminus G=(A_1\cap A_2)^\sim\setminus
E.$$ Inside $E$, we have
\begin{align*}
(\widetilde{A}_1\cap\widetilde{A}_2)\cap E&=\P(N_FA_1)\cap \P(N_FA_2)=\P(T_{A_1}/T_F)\cap\P(T_{A_2}/T_F)\\
&=\P((T_{A_1}\cap T_{A_2})/T_F)=\P(T_{A_1\cap A_2}/T_F)=\P(N_F(A_1\cap A_2))\\
&=(A_1\cap A_2)^\sim\cap E,
\end{align*}
where $N_F(A_1\cap A_2)$ stands for the normal bundle of $F$ in $A_1\cap A_2$. Note that in the fourth equality we have used the condition that $A_1$ and $A_2$ intersect cleanly.

Hence $\widetilde{A}_1\cap\widetilde{A}_2=(A_1\cap A_2)^\sim$.

According to Lemma \ref{ideal_clean}, $\widetilde{A}_1$ and $\widetilde{A}_2$
intersect cleanly if and only if
\begin{equation}\label{e4}
  \I_{\widetilde{A}_1}+\I_{\widetilde{A}_2}=\I_{(A_1\cap A_2)^\sim}.
\end{equation}
But $\widetilde{A}_1=\mathcal{R}(E, \pi^{-1}(A_1))$, where $\mathcal{R}(E, \pi^{-1}(A_1))$ is the residue scheme
to $E$ in $\pi^{-1}(A_1)$ (see \cite{Keel} Theorem 1 or \cite{Fulton}
\S9.2). By a property of residue schemes, we have
$$\I_{\mathcal{R}(E, \pi^{-1}(A_1))}\cdot\I_E=\I_{\pi^{-1}(A_1)},$$ which is the same as
$$\I_{\widetilde{A}_1}\cdot\I_E=\I_{\pi^{-1}(A_1)}.$$ Similarly, we have
$$\I_{\widetilde{A}_2}\cdot\I_E=\I_{\pi^{-1}(A_2)},$$
$$\I_{(A_1\cap A_2)^\sim}\cdot\I_E=\I_{\pi^{-1}(A_1\cap A_2)}.$$
Since $A_1$ and $A_2$ intersect cleanly, $\I_{A_1}+\I_{A_2}=\I_{A_1\cap A_2}$, which implies
$$\I_{\pi^{-1}(A_1)}+\I_{\pi^{-1}(A_2)}=\I_{\pi^{-1}(A_1\cap A_2)}.$$
Thus we get an equality
$$\I_{\widetilde{A}_1}\cdot\I_E+\I_{\widetilde{A}_2}\cdot\I_E=\I_{(A_1\cap
A_2)^\sim}\cdot\I_E.$$ Since $\I_E$ is an invertible sheaf, the
above equality implies (\ref{e4}), hence (iii) is proved.

\end{proof}

\medskip
\begin{proof}[Proof of Proposition \ref{arrangment in blow-up}]
(i) We need to check that any two elements in $\widetilde{\mathcal{S}}$
intersect cleanly, and the intersection is still in $\widetilde{\mathcal{S}}$. For
this, we need to check three cases: ${\widetilde{S}}\cap
{\widetilde{S'}}$,
 ${\widetilde{S}}\cap ({\widetilde{S'}}\cap E)$ and  $({\widetilde{S}}\cap
E)\cap ({\widetilde{S'}}\cap E)$. We only
show the proof for the first case since the proofs for the others are similar.

Suppose $S, S'\in \mathcal{S}$. We can
assume $S\cap S'\neq\emptyset$, otherwise ${\widetilde{S}}\cap {\widetilde{S'}}$
is obviously empty. Suppose the $F$-factorizations of $S$ and $S'$
are $S=A\cap B$ and $S'=A'\cap B'$, respectively. By Lemma \ref{fact
1} (iv),  the $F$-factorization of $S\cap S'$ is $(A\cap A')\cap(B\cap B')$.
Lemma \ref{basic}(v) asserts that ${\widetilde{S}}={\widetilde{A}}\cap {\widetilde{B}}$ and
${\widetilde{S}'}={\widetilde{A}'}\cap {\widetilde{B}'}$. To show that ${\widetilde{S}}$
and ${\widetilde{S'}}$ intersect cleanly along a subvariety in
$\widetilde{\mathcal{S}}$, we consider three cases:

a) $F\subsetneq A\cap A'$. In this case $(S\cap
S')^\sim=(A\cap A')^\sim\cap (B\cap B')^\sim$ and
$${\widetilde{S}}\cap {\widetilde{S'}}=({\widetilde{A}}\cap{\widetilde{A}'})\cap
({\widetilde{B}}\cap{\widetilde{B}'})=(A\cap A')^\sim\cap (B\cap B')^\sim=(S\cap
S')^\sim.$$ Moreover,
$$T_{{\widetilde{S}}}\cap T_{{\widetilde{S'}}}=T_{\widetilde{A}}\cap T_{\widetilde{B}}\cap T_{\widetilde{A}'}\cap
T_{\widetilde{B}'}=T_{(A\cap A')^\sim}\cap T_{(B\cap B')^\sim}=T_{(S\cap
S')^\sim},$$ where the first and third equalities hold because of
Lemma \ref{basic} (v) and the second equality holds because of Lemma
\ref{basic} (iii) and (iv). Thus ${\widetilde{S}}$ intersects ${\widetilde{S'}}$
cleanly along $(S\cap S')^\sim\in \widetilde{\mathcal{S}}$.

b) $F=A\cap A'$ but $F\neq A$ and $F\neq A'$. By Lemma \ref{basic}
(ii), ${\widetilde{A}}\cap{\widetilde{A}'}=\emptyset$, hence
$${\widetilde{S}}\cap
{\widetilde{S}}=({\widetilde{A}}\cap{\widetilde{A}'})\cap
({\widetilde{B}}\cap{\widetilde{B}'})=\emptyset.$$

c) $F=A$ or $A'$. Without loss of generality, we assume $F=A$. Then
$${\widetilde{S}}\cap {\widetilde{S'}}=({\widetilde{A}}\cap{\widetilde{A}'})\cap
({\widetilde{B}}\cap{\widetilde{B}'})=E\cap({\widetilde{A}'}\cap {\widetilde{B}}\cap
{\widetilde{B}'})=E\cap (A'\cap B\cap B')^\sim.$$
By Lemma \ref{basic} (i) and (v),
$$T_{{\widetilde{S}}}\cap T_{{\widetilde{S'}}}=(T_E\cap T_{{\widetilde{B}}})\cap (T_{{\widetilde{A}'}}\cap
T_{{\widetilde{B}'}})=(T_E\cap T_{\widetilde{A}'})\cap T_{(B\cap
B')^\sim}=T_{E\cap (A'\cap B\cap B')^\sim},$$ so again ${\widetilde{S}}$
intersects ${\widetilde{S'}}$ cleanly along $E\cap (A\cap B\cap B')^\sim\in
\widetilde{\mathcal{S}}$. Therefore we have shown that ${\widetilde{S}}$
and ${\widetilde{S'}}$ intersect cleanly along a subvariety in
$\widetilde{\mathcal{S}}$ in all possible cases.

\medskip

\nid(ii) To show that $\widetilde{\mathcal{G}}:=\{\widetilde{G}\}_{G\in\mathcal{G}}$ forms
a building set, we need to show that $\forall{\widetilde{S}}$ (resp. $({\widetilde{S}}\cap E)$) $\in\widetilde{\mathcal{S}}$, the $\widetilde{\mathcal{G}}$-factors of ${\widetilde{S}}$ (resp.
of $({\widetilde{S}}\cap E)$) intersect transversally along ${\widetilde{S}}$
(resp. along $({\widetilde{S}}\cap E)$ ).

By Definition-Lemma \ref{fact 1}, we can assume $S=(G_1\pitchfork
\cdots\pitchfork G_m)\pitchfork(G_{m+1}\pitchfork\cdots\pitchfork
G_k)$, $F\subseteq G_1,\dots, G_m$, and $F\pitchfork
G_{m+1},\dots, G_k$. Define $A=G_1\pitchfork
\cdots\pitchfork G_m$ and $B=G_{m+1}\pitchfork\cdots\pitchfork
G_k$. Then ${\widetilde{S}}={\widetilde{A}}\cap{\widetilde{B}}$ by Lemma
\ref{basic} (v).

Two cases need to be considered: $F\subsetneq A$ and $F=A$. We only give the proof for the first case, since the second case can be proved analogously. Assume $F\subsetneq A$.

First, we show that $\forall{\widetilde{S}}\in\widetilde{\mathcal{S}}$, the $\widetilde{\mathcal{G}}$-factors of ${\widetilde{S}}$ (resp.
of $({\widetilde{S}}\cap E)$) intersect transversally along ${\widetilde{S}}$. Lemma \ref{basic} implies that
$${\widetilde{S}}=\widetilde{G}_1\pitchfork\cdots\pitchfork \widetilde{G}_k$$
and that $\widetilde{G}_1,\dots,\widetilde{G}_k$ are all the
$\widetilde{\mathcal{G}}$-factors of ${\widetilde{S}}$. (Indeed, if some
$\widetilde{G}\in\widetilde{\mathcal{G}}$ contains ${\widetilde{S}}$,
then $G=\pi(\widetilde{G})$ contains $S=\pi({\widetilde{S}})$.
Since $G_1,\dots, G_k$ are all the minimal elements in $\mathcal{G}$
that contain $S$, $G$ contains $G_r$ for some $1\le r\le
k$. The inclusion of their dominant transforms still holds:
$\widetilde{G}\supseteq \widetilde{G}_r$.) Therefore the $\widetilde{\mathcal{G}}$-factors
of ${\widetilde{S}}$ intersect transversally.

Next, we show that $\forall ({\widetilde{S}}\cap E)\in\widetilde{\mathcal{S}}$, the $\widetilde{\mathcal{G}}$-factors of $({\widetilde{S}}\cap E)$ intersect transversally along $({\widetilde{S}}\cap E)$. Noticing that 
$${\widetilde{S}}\cap
E=E\pitchfork{\widetilde{A}}\pitchfork{\widetilde{B}} =E\pitchfork
\widetilde{G}_1\pitchfork\cdots\pitchfork \widetilde{G}_k,$$ 
we assert that $E,
\widetilde{G}_1,\dots,\widetilde{G}_k$ are all the $\widetilde{\mathcal{G}}$-factors of
$({\widetilde{S}}\cap E)$ and the conclusion follows. Indeed, it is enough to show that given any
$\widetilde{G}\in \widetilde{\mathcal{G}}$ containing $({\widetilde{S}}\cap E)$, either $\widetilde{G}=E$ or $\widetilde{G}\supseteq
\widetilde{G}_r$ for some $1\le r\le k$. 
The inclusion $\widetilde{G}\supseteq ({\widetilde{S}}\cap E)$ implies
$G\supseteq (S\cap F)$ by taking the image of $\pi$. By Lemma
\ref{fact 2} (ii), we know that $F$, $G_{m+1}, \dots, G_k$ are all
the $\mathcal{G}$-factors of $(S\cap F)$. Therefore $G$ contains
either $F$ or one of $G_r$ for $m+1\le r\le k$. In the latter case,
we immediate get the conclusion. So we assume that $G$ contains
$F$. If $G=F$, then $\widetilde{G}=E$ from which the conclusion follows. Now we assume $G\supsetneq F$. Since
$$\widetilde{G}\cap E\cong \P(N_FG),
\, \widetilde{S}\cap E\cong\P(N_FA|_{F\cap B}),
$$
and that $\widetilde{G}\cap E$ contains $\widetilde{S}\cap E$, 
we have $(N_FG)_y\supseteq (N_FA)_y$ for any $y\in
F\cap B$. But
$(N_FG)_y\cong T_{G,y}/T_{F,y}$ and $(N_FA)_y\cong
T_{A,y}/T_{F,y}$, therefore $T_{G,y}\supseteq T_{A,y}$ and $G$ contains $A$. Since $G_1,\dots,G_m$
are the $\mathcal{G}$-factors of $A$ by Lemma \ref{fact 2} (i), $G$ contains $G_r$ for some $1\le r\le m$.

\medskip

\nid(iii)``$\mathcal{T}$ is a nest $\Rightarrow$ $\widetilde{\mathcal{T}}$ is a
nest''. Suppose $\mathcal{T}$ is induced by the flag $S_1\subseteq
S_2\subseteq\dots\subseteq
  S_k$. If $S_1\nsubseteq F$ or $S_k\subseteq F$, then $\widetilde{\mathcal{T}}$ is induced by the flag $\widetilde{S}_1\subseteq \widetilde{S}_2
  \subseteq\dots\subseteq
  \widetilde{S}_k$; otherwise there exists an integer $m$, $1\le m\le k-1$ where $S_m\subseteq F$ but $S_{m+1}\nsubseteq F$. In this case
it can be easily checked that $\widetilde{\mathcal{T}}$ is generated by the flag
  \begin{equation*}\label{e8}
  (\widetilde{S}_1\cap \widetilde{S}_{m+1})\subseteq\dots\subseteq
  (\widetilde{S}_m\cap \widetilde{S}_{m+1})\subseteq(\widetilde{S}_{m+1}\cap E)\subseteq\dots\subseteq (\widetilde{S}_k\cap
  E).
  \end{equation*}

  ``$\mathcal{T}$ is a nest $\Leftarrow$ $\widetilde{\mathcal{T}}$ is a nest''. Suppose $\widetilde{\mathcal{T}}$ is induced by the flag $S'_1\subseteq S'_2\subseteq\dots\subseteq
  S'_k$. If $S'_1\nsubseteq E$, then $\mathcal{T}$ is induced by the flag $\pi(S'_1)\subseteq \pi(S'_2)\subseteq
  \dots\subseteq \pi(S'_k)$ and we are done. Now assume $S'_1\subseteq E$, and denote by $m$ the maximal integer satisfying $S'_m\subseteq E$.
  Since $E$ is both minimal and maximal in $\widetilde{\mathcal{G}}$, the $E$-factorization
  of $S'_i$ must be of the form $E\cap B'_i$ for $1\le i\le m$. Then it can be checked that $\mathcal{T}$ is induced by the following flag
$$  (G\cap\pi(B'_1))\subseteq\pi(B'_1)\subseteq\dots\subseteq
  \pi(B'_m)\subseteq\pi(S'_{m+1})\subseteq\dots\subseteq \pi(S'_k).
$$\end{proof}

\begin{proof}[Proof of Proposition \ref{G+}] The proof is similar to the proof of Proposition
\ref{arrangment in blow-up}. The only new case is when $F$ is not minimal in $\mathcal{G}_+$. So throughout the proof we assume that $G_0$ is minimal and $G_0\subsetneq F$. 

\nid(1) We show $\widetilde{\mathcal{S}}_+$ is an arrangement.

First we prove that the intersection $(\widetilde{G}_0\cap\widetilde{S})\cap {\widetilde{S'}}$ is clean for
$S, S'\in\mathcal{S}$. Take the $F$-factorization $S=A\cap B$ and
$S'=A'\cap B'$ in the arrangment $\mathcal{S}$. Take the $G_0$-factorization $B=B_1\cap B_2$ in the arrangement $\mathcal{S}_+$. Similar to the proof of Proposition
\ref{arrangment in blow-up} (i), we need to consider three cases:

a) $F\subsetneq A\cap A'$. Then
\begin{align*}
(\widetilde{G}_0\cap \widetilde{S})\cap {\widetilde{S'}}&=\widetilde{G}_0\cap \widetilde{A} \cap \widetilde{B}_1\cap \widetilde{B}_2\cap \widetilde{A}'\cap \widetilde{B}'\\
&=\widetilde{G}_0\cap \widetilde{A}\cap \widetilde{A}'\cap \widetilde{B}_2\cap \widetilde{B}'\\
&=\widetilde{G}_0\cap (A\cap A')^\sim\cap (B_2\cap B')^\sim\\
&=\widetilde{G}_0\cap (A\cap A'\cap B_2\cap B')^\sim
\end{align*}
The second equality holds because
$\I_{\widetilde{B}_1}=\phi^{-1}\I_{B_1}\subseteq\phi^{-1}
\I_{G_0}=\I_{\widetilde{G}_0}$ therefore $\widetilde{B}_1\supseteq \widetilde{G}_0$.
The third and fourth equalities are because of Lemma \ref{basic}.

Moreover, we have
\begin{align*}
T_{\widetilde{G}_0\cap \widetilde{S}\cap {\widetilde{S}'}}&=T_{\widetilde{G}_0}\cap T_{\widetilde{A}\cap \widetilde{A}'}\cap T_{\widetilde{B}_2\cap {\widetilde{B}'}}\\
&=T_{\widetilde{G}_0}\cap T_{\widetilde{A}}\cap T_{\widetilde{A}'}\cap T_{\widetilde{B}_2}\cap T_{\widetilde{B}'}\\
&=(T_{\widetilde{G}_0}\cap T_{\widetilde{A}}\cap T_{\widetilde{B}_1}\cap T_{\widetilde{B}_2})\cap (T_{\widetilde{A}'}\cap T_{\widetilde{B}'})\\
&=T_{\widetilde{G}_0\cap \widetilde{S}}\cap T_{{\widetilde{S'}}}.
\end{align*}

b) $F=A\cap A'$ but $F\neq A$ and $F\neq A'$. It is easy to
verify that $(\widetilde{G}_0\cap \widetilde{S})\cap {\widetilde{S'}}=\emptyset$.

c) $F=A'$. The proof is similar to a) and we omit it.

Similarly, we can check that $(\widetilde{G}_0\cap \widetilde{S})\cap
(\widetilde{G}_0\cap {\widetilde{S'}})$ and $(\widetilde{G}_0\cap \widetilde{S})\cap (E\cap
{\widetilde{S'}})$ are clean intersections along some elements in
$\widetilde{\mathcal{S}}_+$.

\nid (2) We show $\widetilde{\mathcal{G}}_+$ is a building set. It is enough to show that
the minimal elements in $\widetilde{\mathcal{G}}_+$ which contain $\widetilde{G}_0\cap \widetilde{S}$ intersect transversally along $\widetilde{G}_0\cap
\widetilde{S}$. Assume the $F$-factorization of $S$ is $A\cap B$ where
$G\subsetneq A$. (If $G=A$, then
 $\widetilde{G}_0\subseteq E$ hence we can replace $S$ by $B$ and keep $\widetilde{G}_0\cap \widetilde{S}$
 unchanged.) Assume the $G_0$-factorization of $B$ is $B_1\cap B_2$.

We claim that the set of $\widetilde{\mathcal{G}}_+$-factors of
$\widetilde{G}_0\cap \widetilde{S}$ is:
$$\mathcal{P}:=\{\widetilde{G}_0\}\cup \{\hbox{ $\widetilde{\mathcal{G}}$-factors of } \widetilde{A}\}\cup\{ \hbox{ $\widetilde{\mathcal{G}}$-factors of }
\widetilde{B}_2\}.$$ It is easy to check that the subvarieties in $\mathcal{P}$ intersect transversally. To
show that the subvarieties in $\mathcal{P}$ are all the $\widetilde{\mathcal{G}}_+$-factors, it suffices to show that any minimal element $\widetilde{G}\in\widetilde{\mathcal{G}}_+$ that
contains $\widetilde{G}_0\cap \widetilde{S}$ ($=\widetilde{G}_0\cap \widetilde{A}\cap
\widetilde{B}_2$) belongs to $\mathcal{P}$.

Since $G=\phi(\widetilde{G})\supseteq\phi(\widetilde{G}_0\cap \widetilde{A}\cap
\widetilde{B}_2)=G_0\cap B_2$ and the $\mathcal{G}_+$-factors of
$G_0\cap B_2$ are $G_0$ and $\mathcal{G}$-factors of $B_2$,
$G\supseteq G_0$ or $G\supseteq B_2$. If $G\supseteq B_2$, then the conclusion follows, so we can assume $G\supseteq G_0$. Then $G\pitchfork F$ or $G\supseteq F$.

If $G\pitchfork F$, then $G\supseteq G_0$ implies 
$\widetilde{G}\supseteq \widetilde{G}_0$. Since $\widetilde{G}$ is chosen to be minimal, $\widetilde{G}=\widetilde{G}_0$ belongs to $\mathcal{P}$.

If $G\supseteq F$, then $\widetilde{G}=\P(N_FG)$ and $\widetilde{G}_0\cap
\widetilde{S}=\P(N_FA|_{G_0\cap B})$. Therefore $G\supseteq A$ which
implies that $G$ is a $\mathcal{G}$-factor of $A$ and $\widetilde{G}$
is a $\widetilde{\mathcal{G}}$-factor of $\widetilde{A}$, hence $\widetilde{G}\in \mathcal{P}$.
\end{proof}

\sss {Codimensions of the centers} (Thanks the referee for suggesting this question.)
In the original construction of the Fulton-MacPherson configuration space $X[n]$, each blow-up is along a nonsingular center of codimension $m$ or $m+1$, where $m$ is the dimension of $X$. But if we construct $X[n]$ by blowing up centers in ascending dimension, the codimensions of the centers are much larger: the first blow-up is along the smallest diagonal $\Delta_{[n]}$ which is of codimension $m(n-1)$.

In general, given a specified order of blow-ups, we can find the dimension, hence the codimension, of the centers easily.

\begin{prop} Let $Y$ be a nonsingular variety and $\mathcal{G}=\{G_1,\cdots,G_N\}$ be a nonempty building
set of subvarieties of $Y$ satisfying the condition (*) in Theorem \ref{order}. Let $j$ be an integer between $1$ and $N$. Define $\mathcal{G}':=\{G_1,\dots, G_{j-1}\}$ and define $\mathcal{F}=\{G_{i_1},G_{i_2},\dots,G_{i_\ell}\}$ to be the minimal elements of $\{G\in\mathcal{G}': G\supseteq G_j\}$. 

Then in the construction of $Y_{\mathcal{G}}$ by blowing up along $G_1,\cdots,G_N$ in order, the center of the $j$-th blow-up is of dimension 
$$\dim G_j+\sum_{k=1}^\ell(d-1-\dim G_{i_k})$$ if $\mathcal{F}\neq\emptyset$, and is of dimension $\dim G_j$ if $\mathcal{F}=\emptyset$. 
\end{prop}
\begin{proof} The set $\mathcal{G}':=\{G_1,\dots, G_{j-1}\}$ is a building set by the condition (*) in Theorem \ref{order}. By the same theorem we can assume that $G_1,\dots,G_{j-1}$ is in order of ascending dimension. 
 Denote the variety obtained after the $i$-th blow-up by $Y_i$. We want to find the dimension of $\widetilde{G}_j$ in $Y_{j-1}$. 

Since (a) blowing up a center that does not contain $\widetilde{G}_j$ will not change the dimension of $\widetilde{G}_j$, and (b) $\widetilde{G}$ does not contain $\widetilde{G}_j$ if $G$ does not contain $G_j$,  we can focus on the subset $\mathcal{G}''\subseteq\mathcal{G}'$ of subvarieties that contain $G_j$. 
Let $\mathcal{F}:=\{G_{i_1},G_{i_2},\dots,G_{i_\ell}\}$ be the set of minimal elements in $\mathcal{G}''$. Define $S:=G_{i_1}\cap G_{i_2}\cap\dots\cap G_{i_\ell}$. Then $\mathcal{F}$ is also the set of minimal elements in $\mathcal{G}'$ that contain $S$. By the definition of building set,  $G_{i_1}\dots G_{i_\ell}$ intersect transversally. A subvariety $G\in \mathcal{G}''\setminus \mathcal{F}$ must contain a subvariety, say $G_i$, in $\mathcal{F}$. Then $G\supsetneq G_i\supsetneq G_j$, and  in the variety $Y_{i-1}$ we have $\widetilde{G}\supsetneq \widetilde{G}_i\supsetneq \widetilde{G}_j$. It can be easily checked that in the variety $Y_j$ (obtained by blowing up along $\widetilde{G}_i$), $\widetilde{G}$ will no longer contain $\widetilde{G}_j$, hence the blow-up along $\widetilde{G}$ will not change the dimension of $\widetilde{G}_j$. In other words, we only need to find the change of $\dim \widetilde{G}_j$ for the blow-ups along the transversal subvarieties in $\mathcal{F}$, which is simply
$$\sum_{k=1}^\ell(d-1-\dim G_{i_k})$$
\end{proof} 

\eg Let $m=\dim X$. In the original construction of the Fulton-MacPherson configuration space $X[4]$ (cf. \S\ref{eg:FM}), the dimension of the first center is $\dim\Delta_{12}=3m$, the dimension of the second senter  $\widetilde{\Delta}_{123}$ is $$\dim\Delta_{123}+(4m-1-\dim\Delta_{12})=2m+(m-1)=3m-1,$$ 
it can be easily checked that $\dim\widetilde{\Delta}_I$ is $3m$ if $|I|=2$ and $3m-1$ otherwise, so the codimension is either $m$ or $m+1$. In general, using the order of blow-ups
in the orginal construction of $X[n]$ for $n\ge 2$, the codimension of each center $\widetilde{\Delta}_I$ is $m$ if $|I|=2$ and $m+1$ othewise.

\eg The codimension of the each blow-up center is 2 in Keel's construction of $\overline{M}_{0,n}$. Since the blow-ups in Keel's construction of $\overline{M}_{0,n}$ can be obtained by restricting the blow-ups in the original Fulton-MacPherson's construction of $X[n]$ to a fiber of $\pi_{123}$ (defined in Corollary \ref{corM0n}), each blow-up center is of codimension $m=1$ or $m+1=2$. But blowing up along a center of codimension 1 does nothing. So we only need to carry out blow-ups along codimension-2 centers.

\sss{The statements for a general arrangement}\label{general}
\begin{df}
A {\em arrangement} of subvarieties of a nonsingular variety\, $Y$ is a finite set $\mathcal{S}=\{S_i\}$ of nonsingular closed
subvarieties $S_i$ properly contained $Y$ satisfying the following conditions:
\begin{itemize}
\item[(i)] $S_i$ and $S_j$ intersect cleanly,
\item[(ii)] $S_i\cap S_j$ either equals to a disjoint union of some $S_k$'s or is empty.
\end{itemize}
\end{df}

\begin{df}
 Let $\mathcal{S}$ be an arrangement of subvarieties of $\, Y$. A subset $\mathcal{G}\subseteq \mathcal{S}$ is
called a {\em building set of $\mathcal{S}$} if there is an open cover $\{U_i\}$ of\, $Y$ such that the restriction of the arrangement $\mathcal{S}|_{U_i}$ is simple  for each $i$ and  $\mathcal{G}|_{U_i}$ is a building set of $\mathcal{S}|_{U_i}$.

A finite set $\mathcal{G}$ of nonsingular subvarieties of $Y$ is
called a {\em building set} if the set of all possible
intersections of collections of subvarieties from $\mathcal{G}$
forms an arrangement $\mathcal{S}$, and that $\mathcal{G}$ is a
building set of $\mathcal{S}$ defined as above. In this situation,
$\mathcal{S}$ is called the {\em arrangement induced by
$\mathcal{G}$}.
\end{df}

\begin{df}
A subset $\mathcal{T}\subseteq\mathcal{G}$ is called
{\em $\mathcal{G}$-nested} (or a {\em $\mathcal{G}$-nest}) if
there is an open cover $\{U_i\}$ of $Y$, such that the restriction of the arrangement induced by $\mathcal{G}$ to each $U_i$ is simple and $\mathcal{T}|_{U_i}$ is a $\mathcal{G}|_{U_i}$-nest.
\end{df}

We define the wonderful compactification same as Definition \ref{def wonderful compactification}. Then
Theorem \ref{main wonder} and Theorem \ref{order} still hold.

\bigskip

\noindent {\sc Li Li }\\
Department of Mathematics\\
University of Illinois at Urbana-Champaign\\
Email: {\tt llpku@math.uiuc.edu}\\

\end{document}